%% file: sQ.tex
\documentclass{article}[11pt]
\usepackage{geometry}
\usepackage{enumerate}
\usepackage{latexsym,booktabs}
\usepackage{amsmath,amssymb, bm}
\usepackage{graphicx}
\usepackage{subfigure}
\usepackage{amsthm}
\usepackage{thmtools}
\usepackage{float}        
\usepackage{longtable}
\usepackage{multirow}
\usepackage{diagbox}
\usepackage{subfigure}
\usepackage{natbib}        
\usepackage{url}
\usepackage{textcomp}      
\usepackage{color}
\usepackage{subfiles}

\renewcommand{\vec}[1]{\boldsymbol{#1}} 

\usepackage{siunitx} %
\usepackage[flushleft, online]{threeparttablex} %

\usepackage{algorithmicx}
\usepackage[ruled]{algorithm}
\usepackage{algpseudocode}

\usepackage{cases}

\renewcommand{\arraystretch}{1.5}

\theoremstyle{definition}
\newtheorem{exmp}{Example}

\newtheorem*{Proof}{Proof}
\newtheorem{Lemma}{Lemma}

\renewcommand\thmcontinues[1]{\textbf{Continued}}

\usepackage{authblk}

\title{MILP Approximations for non-stationary stochastic lot-sizing under ($s,Q$)-type policy}

\author[a]{Xiyuan Ma \thanks{Corresponding author: Xiyuan.Ma@ed.ac.uk}}
\author[a]{Roberto Rossi}
\author[a]{Thomas Welsh Archibald}
\affil[a]{Business School, University of Edinburgh, Edinburgh, United Kingdom}

\date{}

\begin{document}
\maketitle

\begin{abstract}
    This paper addresses the single-item single-stocking location non-stationary stochastic lot-sizing problem under a reorder point -- order quantity control strategy. The reorder points and order quantities are chosen at the beginning of the planning horizon. The reorder points are allowed to vary with time and we consider order quantities either to be a series of time-dependent constants or a fixed value; this leads to two variants of the policy: the ($s_t,Q_t$) and the ($s_t,Q$) policies, respectively. For both policies, we present stochastic dynamic programs (SDP) to determine optimal policy parameters and introduce mixed integer non-linear programming (MINLP) heuristics that leverage piecewise-linear approximations of the cost function. Numerical experiments demonstrate that our solution method efficiently computes near-optimal parameters for a broad class of problem instances.\\

\end{abstract}

\section{Introduction}
The non-stationary stochastic lot-sing problem is an extension of the well-known dynamic lot-sizing problem \citep{wagner1958dynamic}. In this problem, one considers a single-item single-stocking location inventory system under a finite planning horizon and periodic review; the demand is stochastic and non-stationary. To deal with the uncertainty inherent in a stochastic lot-sizing problem, \citet{bookbinder1988strategies} introduced three control strategies: the ``static uncertainty,'' the ``static-dynamic uncertainty,'' and the ``dynamic uncertainty,'' which represent different approaches for determining the timing and size of orders.

\citeauthor{bookbinder1988strategies}'s control strategies are captured by various policies. The ($R,Q$) policy determines the inventory review schedule $R$ and the order quantity $Q$ before the system operates; this is the static uncertainty strategy. The ($s,S$) policy is the dynamic uncertainty strategy, in which the timing and size of orders are decided as late as possible, in a wait-and-see fashion, by leveraging the reorder point $s$, and the order-up-to level $S$. \citet{scarf1959optimality} showed that if the holding and shortage costs are convex, the optimal policy in each period is of ($s,S$) type. In a static-dynamic uncertainty strategy one either fixes at the set the order schedule, and computes the exact order quantity only when orders are issued, via suitable order-up-to-levels; or fixes the order quantities at the set, and decides when orders are issued in a wait-and-see fashion, by relying on a reorder threshold. This leads to the ($R,S$) policy and ($s,Q$) policy (also referred to as the ($r,Q$) policy), respectively.

Compared to stationary demand, there are relatively few studies in the literature that consider non-stationary demand. However, in the majority of practical circumstances, demand is not only stochastic but also non-stationary.

In research on the ($R,Q$) policy for non-stationary demand, \citet{sox1997dynamic} proposes a MINLP of the dynamic lot-sizing problem with dynamic costs and develops a solution algorithm that resembles the Wagner-Whitin algorithm. This policy is also investigated by \citet{vargas2009optimal}, who develops a stochastic dynamic programming model which is equivalent to a shortest path problem in a specified acyclic network. \citeauthor{vargas2009optimal} also provides an optimisation algorithm with rolling horizon with two stages: (1) to determine optimal replenishment quantities for any sequence of replenishment points, and (2) to identify the optimal sequence of replenishment points.

For the static-dynamic uncertainty strategy, research under non-stationary demand mostly considers the ($R,S$) policy. \citet{tarim2004stochastic} formulates the problem as a mixed integer program (MIP). They model the total expected cost by minimising the summation of holding and ordering costs under a constraint on the probability of the closing inventory in each time period. A method to solve this model efficiently is introduced in \citep{tarim2011efficient}, where the relaxation of the original MIP model is converted to a shortest path problem and implemented by branch-and-bound procedures. \citet{tarim2006modelling} provide another MIP formulation where the objective function is obtained by the mean of a piecewise linearisation. The accuracy of the approximation can be adjusted ad libitum by introducing new breakpoints.

\citet{ozen2012static} consider both penalty cost and service level and prove that the optimal policy is a base stock policy for both penalty and service-level constrained models, and also for capacity limitations and minimum order quantity requirements. More recently, \citet{rossi2015piecewise} consider several service level measures --- $\alpha$ service level on each period, $\beta^{cyc}$ service level independently for each replenishment cycle, and the classic $\beta$ service level --- by adding suitable constraints that leverage the loss function and its complementary function to describe the expected total holding and penalty cost. A piecewise linearisation approach is utilized to convert the cost function from non-linear to linear form.

Computing ($s,S$) policy parameters under non-stationary demand is a challenging task. The classic Silver and Meal heuristic algorithm \citep{silver1973heuristic} for deterministic demand has been extended by \citet{silver1978inventory} and \citet{askin1981procedure}. \citeauthor{silver1978inventory}'s algorithm uses a deterministic model to calculate the number of periods that each order must cover; when this replenishment plan is known, the associated safety stocks are then myopically determined. \citet{askin1981procedure} explicitly includes the cost effects of probabilistic demand in the choice of the number of periods in which to order.
\citet{bollapragada1999simple} approximate the non-stationary problem via a series of stationary problems based on the method developed by \citet{zheng1991finding}. Parameters are determined by equating the cumulative mean demand of stationary and non-stationary problems over the expected reorder cycle. \citet{xiang2018computing} introduce a MINLP formulation for $(s,S)$ policy by applying the piecewise linearisation approximation proposed by \citet{rossi2015piecewise}. \citeauthor{xiang2018computing} also derive a heuristic algorithm with binary search. Both solution methods outperform the previous heuristics in computational efficiency for short and long time horizon tests. The comparison between two proposed algorithms shows that binary search requires significantly less time than the MINLP.

Based on this literature survey, we note a gap in the study of non-stationary demand: no literature discussed or investigated the static-dynamic uncertainty strategy in the form of an ($s,Q$) policy. In this paper, we focus on the stochastic lot-sizing problem under non-stationary demand and an ($s,Q$) control strategy. The reorder points $s_t$ vary with time, and we consider two cases of order quantity, which is either able to shift according to the time periods ($Q_t$) or is fixed over the planning horizon ($Q$). This leads to two ($s,Q$)-type policies: the ($s_t,Q_t$) policy and the ($s_t,Q$) policy. These policies involve determining $s_t$ and $Q_t$ (or $Q$) values at the beginning of the planning horizon.

Compared to the optimal policy introduced by \citet{scarf1959optimality} which allows the order quantity to vary with inventory level and time period, the order quantity in an ($s_t,Q_t$) policy is only affected by the time period and applies to all inventory levels, while the order quantity in an ($s_t,Q$) policy is a constant value for the entire planning horizon, and does not shift with inventory level or time period.

We make the following contributions to the stochastic lot-sizing literature.
\begin{itemize}\setlength{\itemsep}{-0.05cm}
\item We model the non-stationary stochastic lot-sizing problem under a static-dynamic uncertainty policy in which order quantities are determined ``statically'', at the onset of the planning horizon, while reordering decisions are determined ``dynamically'', in a wait-and-see-fashion.{\em We prove that the resulting optimal policy takes the non-stationary} ($s,Q$) {\em form}.
\item To efficiently determine near-optimal policy parameters of the non-stationary ($s_t,Q_t$) and ($s_t,Q$) policies, we present a heuristic algorithm  based on a MINLP and binary-search. The model is then turned into a mixed integer linear program by applying the piecewise linearisation approach discussed in \citep{rossi2014piecewise}.
\item In a comprehensive numerical study, based on instances drawn from \citep{xiang2018computing}, we investigate the performance of the ($s_t,Q_t$) and ($s_t,Q$) policies against an optimal ($s_t,S_t$) policy. We show that optimality gaps for the ($s_t, Q_t$) policy obtained via our heuristic are tighter than those of a near-optimal ($R_t, S_t$) policy obtained via the approach in \citet{rossi2015piecewise}. Finally, we observe that an ($s_t, Q$) policy lacks flexibility and leads to substantial optimality gaps.
\end{itemize}

The rest of this paper is structured as follows. In Section \ref{sec:problemDescription} we introduce the problem settings and present a stochastic dynamic programming (SDP) formulation. Section \ref{sec:sdps} discusses the stochastic dynamic programming formulation of the ($s_t, Q_t$) and ($s_t,Q$) policies. We also show that the resulting optimal policies take the non-stationary ($s_t,Q_t$) and ($s_t,Q$) forms through the uniqueness of reorder points. In Section \ref{sec:heuristic}, we apply an existing MINLP model and a binary search approach to the ($s,S$) policy, based on which we derive a heuristic algorithm to compute near-optimal policy parameters of the ($s_t,Q_t$) policy and discuss the application of this algorithm on the ($s_t, Q$) policy. A computational analysis is presented in Section \ref{sec:computations} and we reach our conclusions in Section \ref{sec:conclusion}.


\section{Problem description}\label{sec:problemDescription}
We consider a single-item single-location non-stationary stochastic lot-sizing problem over a planning horizon of $T$ periods. Replenishment orders are placed and instantaneously delivered at the beginning of each time period. Each replenishment order incurs an ordering cost $c(\cdot)$ comprising a fixed ordering cost $K$ and a linear ordering cost $z$ proportional to the non-negative order quantity $Q$, where
\begin{equation}\label{eq:OrderingCost}
\renewcommand\arraystretch{1}
c(Q)\triangleq
\left\{
  \begin{array}{ll}
    \hbox{$K+z\cdot Q$}, & \hbox{$Q>0$;} \\
    0, & \hbox{$Q=0$.}
  \end{array}
\right.
\end{equation}
The periods' demands $d_t$, for $t=1,\cdots,T$, are independent random variables with known probability density functions $g_t(\cdot)$. Any unmet demand at the end of the period is back-ordered. At the end of each period, a linear holding cost $h$ is incurred for each unit carried from one period to the next, and a linear penalty cost $b$ is charged on each unit back-ordered. The expected immediate holding and penalty cost at the end of period $t$ is expressed as
\begin{equation}\label{eq:expectedImmediateCost}
L_t(y) \triangleq \mathbb{E}[h\max(y-d_t) + b\max(d_t-y)],
\end{equation}
where $y$ denotes the inventory level after receiving the replenishment and $\mathbb{E}[\cdot]$ denotes the expectation operator.

Let $C_t(x)$ represent the expected total cost of an optimal policy over periods $t,\ldots,T$ with opening inventory level $x$; then the problem can be modelled as a stochastic dynamic program \citep{BellmanRichard1957Dp}
\begin{equation}\label{eq:sS-sdp-objFunc}
C_t(x) \triangleq \min_{y\geq x}\mathop{}\{c(y-x) + L_t(y) +\mathbb{E}[C_{t+1}(y - d_t)]\},
\end{equation}
where $C_{T+1}(x) \triangleq 0$, is the boundary condition.

\citet{scarf1959optimality} showed that, if $L_t(y)$ is convex, the optimal policy of the dynamic inventory problem is of an ($s,S$) type, where the inventory system places a replenishment to reach the order-up-to level $S$ when the stock is found to be below the reorder point at a review point. This conclusion is based on a study of the function $G_t(y) + zy$, where
\begin{equation}\label{eq:sdp-optimal-noOrder}
G_t(y) \triangleq L_t(y) + \mathbb{E}[C_{t+1}(y-d_t)],
\end{equation}
and $G_t(y)$ represents the expected total cost over period $t$ to $T$ when the opening inventory is $y$ and no order is placed in period $t$. Table \ref{tab:App.notations} in Appendix \ref{sec:App.notations} summarises the notation functions used in this paper.

In the rest of this paper, we conduct the discussion assuming $L_t(y)$ convex. In fact, as the holding and penalty costs used in this paper are linear, $L_t(y)$ is a weighted sum of two convex functions and hence convex. A detailed proof can be found in \citep[page~490]{rossi2014piecewise}.

\begin{exmp}[label=ex1]\label{ex:classicLarge}
Consider a 4--period stochastic lot-sizing problem under Poisson-distributed demand with rates  $d_t = \langle20,40,60,40\rangle$ . The cost parameters are $K=100$, $z=0$, $h=1$ and $b=10$. Fig. \ref{fig:ex1.sS} illustrates the variation of $G_t(I_0)$ with $I_0\in [0,200]$ and no replenishment order placed in period 1, where $G_1(0) = 481$.

\begin{figure}[htb]
\centering
\includegraphics[width=0.8\textwidth]{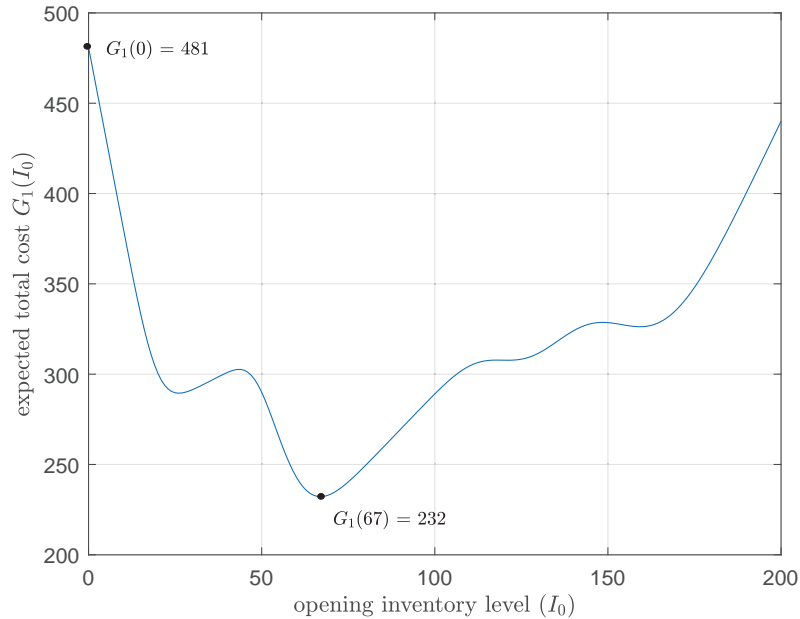}
\caption{Plot of $G_1(I_0)$}\label{fig:ex1.sS}
\end{figure}
\end{exmp}


\section{Stochastic dynamic programs of ($\bm{s_t}$,$\bm{Q_t}$) and ($\bm{s_t}$,$\bm{Q}$) policy}\label{sec:sdps}
This section introduces the stochastic dynamic programming formulations of the stochastic lot-sizing problem under the ($s_t,Q_t$) policy and the ($s_t,Q$) policy in Section \ref{sec:sdp-sQt} and Section \ref{sec:sdp-sQ-Q}, respectively.

\subsection{A stochastic dynamic program for (\textit{$\bm{s_t,Q_t}$}) policy}\label{sec:sdp-sQt}
An $(s_t,Q_t)$ policy places a replenishment order of size $Q_t$ at the beginning of period $t$ if the inventory level is below the reorder point $s_t$, and does not place any order otherwise \citep{silver1998inventory}. The optimal expected total cost of the system controlled under an ($s_t,Q_t$) policy can be determined by computing all feasible combinations of reorder quantity $Q_t$, for $t=1,\ldots, T$. Let $\vec{q}_t = \langle Q_t, \ldots, Q_T\rangle$ denote a $(T-t+1)$-dimensional vector representing order quantities $Q_t, \ldots, Q_T$ and $\mathcal{Q}_t$ be the vector space representing all combinations of order quantities $\vec{q}_t$. For any $\vec{q}_t\in \mathcal{Q}_t$, the expected total cost when the opening inventory level is $x$ is denoted as
\begin{equation}
V_t(x,\vec{q}_t) \triangleq \min_{\delta\in\{0,1\}} \mathop{} \{ c(\delta Q_t) + L_t(x+\delta Q_t) +\mathbb{E}[V_{t+1}(x+\delta Q_t - d_t, \vec{q}_{t+1})]\},\label{eq:sQt-sdp-singleComb-t}
\end{equation}
where $\delta$ is a binary variable that represents the reordering decision in period $t$ when the initial inventory level is $x$; finally,
\begin{equation}
V_T(x,\vec{q}_T) \triangleq \min_{\delta\in\{0,1\}} \mathop{} \{ c(\delta Q_T) + L_T(x+\delta Q_T)\}\label{eq:sQt-sdp-singleComb-T}
\end{equation}
is the boundary condition. Therefore, considering all combinations, the optimal expected total cost when the initial inventory level at the beginning of the planning horizon is $x$ can be defined as
\begin{equation}
V_0(x) \triangleq \min_{\vec{q}_1 \in \mathcal{Q}_1} \{ V_1(x, \vec{q}_1)\}.\label{eq:sQt-sdp-optimal}
\end{equation}
Let the optimal order quantity be represented by the vector $\vec{q}_t^* \triangleq \langle Q_t^*,\ldots,Q_T^* \rangle$.

Next we show that the policy found by the formulation in Section \ref{sec:sdp-sQt} is of an ($s_t,Q_t$) form. The following discussion is inspired by the work of \cite{gallego2004all} on all-or-nothing ordering policies under a capacity constraint.
For any opening inventory level $x$ and a vector of order quantities $\vec{q}_t$, let $J_t(x,\vec{q}_t)$ and $\hat{J}_t(x,\vec{q}_t)$ denote the expected total cost when the decision in period $t$ is not to order ($\delta=0$) and to order ($\delta = 1$) respectively, it follows that
\begin{equation}\label{eq:sQt-optimal-case-noReorder}
J_t(x,\vec{q}_t) \triangleq L_t(x) + \mathbb{E}[V_{t+1}(x - d_t, \vec{q}_{t+1})]
\end{equation}
and
\begin{equation}\label{eq:sQt-optimal-case-reorder}
\hat{J}_t(x,\vec{q}_t) \triangleq c(Q_t) + L_t(x+Q_t) + \mathbb{E}[V_{t+1}(x +Q_t - d_t, \vec{q}_{t+1})].
\end{equation}
Recall that Eq.\eqref{eq:sQt-sdp-singleComb-t} optimises the system over the reorder decision $\delta\in\{0,1\}$ and is equivalent to
\begin{eqnarray}
V_t(x,\vec{q}_t) & = & \min\{\hat{J}_t(x,\vec{q}_t),\mathop{} J_t(x,\vec{q}_t)  \}\notag\\
&=& \min\{K + zQ_t+L_t(x+Q_t) + \mathbb{E}[V_{t+1}(x+Q_t-d_t,\vec{q}_{t+1})],\notag\\
&&\hspace{2.5cm} L_t(x) + \mathbb{E}[V_{t+1}(x-d_t,\vec{q}_{t+1})]\}\notag\\
&=&\min\{K+zQ_t+J_t(x+Q_t,\vec{q}_t),\mathop{} J_t(x,\vec{q}_t)\}\notag\\
&=&J_t(x,\vec{q}_t) + \min\{K + zQ_t-\Delta J_t(x,\vec{q}_t),\mathop{} 0\},\label{eq:sQt-s-reviseV}
\end{eqnarray}
where we define
\begin{equation}
\Delta J_t(x,\vec{q}_t) \triangleq J_t(x,\vec{q}_t)- J_t(x+Q_t,\vec{q}_t).
\end{equation}
From Eq.\eqref{eq:sQt-s-reviseV}, it is optimal to reorder in period $t$ with opening inventory $x$ when $\Delta J_t(x,\vec{q}_t)>K+zQ_t$ and not to reorder otherwise. If we choose not to reorder when $\Delta J_t(x,\vec{q}_t) = K+zQ_t$, then the region of opening inventory level $x$ that is optimal to reorder can be expressed as
\begin{equation}\label{eq:regionCompare}
\{x:\Delta J_t(x,\vec{q}_t) > K+zQ_t\}.
\end{equation}
If $\Delta J_t(x,\vec{q}_t)$ is non-increasing in $x$ for an order quantities $\vec{q}_t^*$, then either there exits an $s_t$ such that it is optimal to order in period $t$ when $x<s_t$ and not otherwise, or it is never optimal to order in period $t$; and it hence leads to the ($s_t, Q_t$) policy. In the following, for any given $\vec{q}_t$, we show the monotonicity of $\Delta J_t(x,\vec{q}_t)$ in $x$.


\begin{Lemma}\label{lemma:Ft-non-increasing}
$L_t(y) - L_t(y+a)$ is non-increasing in $y$ for any $a>0$ and $t = 1,\ldots,T$.
\end{Lemma}
\begin{Proof}
Since that $L_t(y)$ is convex, then its derivative $L_t^{'}(y)$ is non-decreasing by the definition of convexity. For any $a>0$ and any $t=1,\ldots,T$,
$[L_t(y) - L_t(y+a)]^{'} = L_t^{'}(y) - L_t^{'}(y+a))\leq 0$; therefore, $L_t(y) - L_t(y+a)$ is non-increasing in $y$. \qedsymbol
\end{Proof}

\begin{Lemma}\label{lemma:Gt-non-increasing}
For a given $\vec{q}_t$, the function $\Delta J_t(x,\vec{q}_t)$ is monotonically non-increasing with respect to the opening inventory level $x$ for any $t=1,\ldots,T$.
\end{Lemma}
\begin{Proof}
We prove this by induction. For period $T$,
\begin{equation}
\Delta J_T(x,\vec{q}_T)  =  J_T(x,\vec{q}_T) - J_T(x+Q_T,\vec{q}_T)
= L_T(x) - L_T(x+Q_T)\notag
\end{equation}
is non-increasing by Lemma \ref{lemma:Ft-non-increasing}. Assuming that $\Delta J_t(x,\vec{q}_t)$ is non-increasing in $x$, we want to show that $\Delta J_{t-1}(x,\vec{q}_{t-1})$ is non-increasing in $x$. We find that
\begin{eqnarray}
&&K+zQ_t+V_t(x+Q_t,\vec{q}_t) - V_t(x,\vec{q}_t)\notag\\
&=&K+zQ_t+J_t(x+Q_t,\vec{q}_t) - J_t(x,\vec{q}_t) + \min\{0,K + zQ_t - \Delta J_t(x+Q_t,\vec{q}_t)\}\notag\\
&& \hspace{5.65cm}- \min\{0,K + zQ_t - \Delta J_t(x,\vec{q}_t)\}\notag\\
&=&K+zQ_t-\Delta J_t(x,\vec{q}_t) + \min\{0,K+zQ_t-\Delta J_t(x+Q_t,\vec{q}_t)\} -\min\{0,K+zQ_t-\Delta J_t(x,\vec{q}_t)\}\notag\\
&=&\max\{0,K+zQ_t-\Delta J_t(x,\vec{q}_t)\} + \min\{0,K + zQ_t - \Delta J_t(x+Q_t,\vec{q}_t)\}\notag
\end{eqnarray}
is the sum of two non-decreasing functions because $\Delta J_t(x,\vec{q}_t)$ is assumed to be non-increasing, then $V_t(x,\vec{q}_t) - V_t(x+Q_t,\vec{q}_t)$ is non-increasing. Consequently, with a non-increasing $L_{t-1}(x) - L_{t-1}(x+Q_{t-1})$ in $x$,
\begin{eqnarray}
\Delta J_{t-1}(x,\vec{q}_t) & = & J_{t-1}(x,\vec{q}_{t-1}) - J_{t-1}(x+Q_{t-1},\vec{q}_{t-1})\notag\\
&=&L_{t-1}(x) - L_{t-1}(x+Q_{t-1}) + \mathbb{E}[V_t(x-d_{t-1},\vec{q}_t) - V_t(x+Q_t-d_{t-1},\vec{q}_t)]\notag
\end{eqnarray}
is the sum of two non-increasing functions; therefore, $\Delta J_{t-1}(x,\vec{q}_t)$ is non-increasing in $x$. This completes the proof by induction. \qedsymbol
\end{Proof}
For a given $\vec{q}_t$, the monotonicity of $\Delta J_t(x,\vec{q}_t)$ in $x$ assures the unique existence of the reorder point $s_t$, which defines the region of opening inventory $x<s_t$ for which it is optimal to reorder, where $s_t$ can be denoted as
\begin{equation}\label{eq:sQt-s}
s_t = \inf\{x:\mathop{} \Delta J_t(x,\vec{q}_t) < K+zQ_t\};
\end{equation}
if the inventory levels are discrete, then $s_t$ is the minimum value of $x$ such that $\Delta J_t(x,\vec{q}_t) < K+zQ_t$, where $Q_t$ is the first argument of the order quantities $\vec{q}_t$. The reorder points associated with the optimal order quantities $\vec{q}_t^*$ hence can be denoted as $\vec{s}_t^*\triangleq \langle s_t^*, \ldots,s_T^* \rangle$.

\begin{exmp}[label=ex2]\label{ex:smallPoisson_stQt}
Consider a 4--period stochastic lot-sizing problem under Poisson-distributed demand with rates $d_t = \langle2, 1, 5, 3\rangle$. The cost parameters are $K= 5$, $z = 0$, $h=1$ and $b = 3$. The maximum order quantity is set to $9$. After exhaustive enumeration of all order quantity vectors, we obtain $\vec{q}_1^* = \langle3, 3, 8, 5\rangle$ and the associated reorder points $\vec{s}_1^*=\langle1, 0, 4, 1\rangle$. The expected total cost of the optimal ($s_t,Q_t$) policy is $22.5$ when the initial inventory is 0.
Under discrete inventory levels with Poisson demand, Fig. \ref{fig:ex2_diff} and \ref{fig:ex2_compare} illustrate determining $s_1^*$ by scatter plots. In Fig. \ref{fig:ex2_diff}, $s^*_1=1$ is selected as the minimum value such that $\Delta J_1(I_0,\vec{q}_1^*) < K$, which is equivalent to $J_1(I_0,\vec{q}^*) > \hat{J}_1(I_0,\vec{q}^*)$ when $I_0\leq 0$, suggesting it is optimal to order; and $J_1(I_0,\vec{q}^*) < \hat{J}_1(I_0,\vec{q}^*)$ when $I_0\geq 1$, suggesting it is optimal not to order, as  Fig. \ref{fig:ex2_compare} shows.
\begin{figure}[H]
    \centering
    \includegraphics[width=0.8\textwidth]{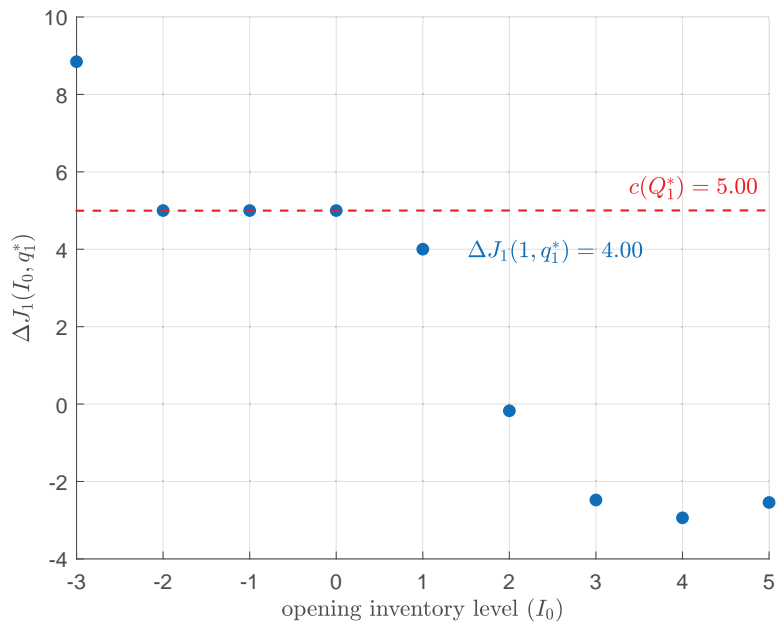}
    \caption{$s^*_1=1$ determined by comparing $\Delta J_1(I_0,\vec{q}_1^*)$ and $c(Q_1^*)$.}
    \label{fig:ex2_diff}
\end{figure}

\begin{figure}[H]
\centering
\includegraphics[width=0.8\textwidth]{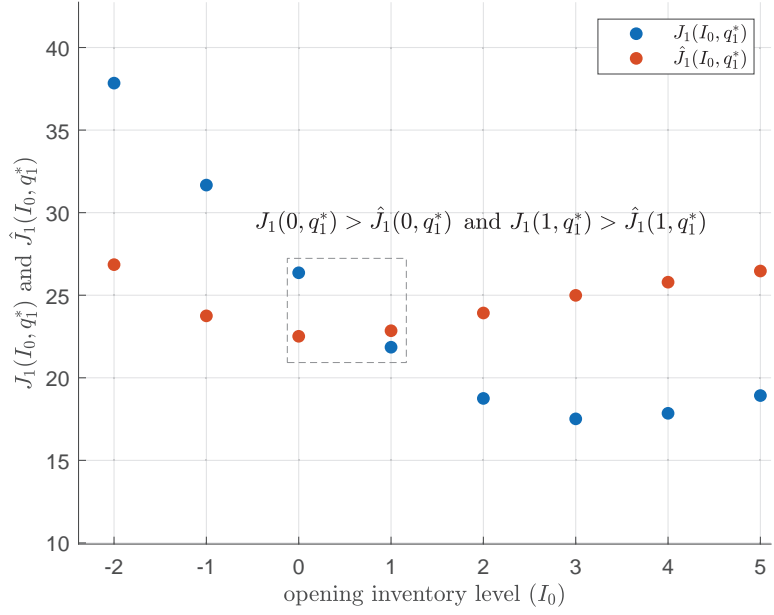}
\caption{$s^*_1=1$ determined by comparing $J_1(I_0,\vec{q}_1^*)$ and $\hat{J}_1(I_0,\vec{q}_1^*)$.}\label{fig:ex2_compare}
\end{figure}

\end{exmp}

\subsection{A stochastic dynamic program for (\textit{$\bm{s_t,Q}$}) policy}\label{sec:sdp-sQ-Q}
An ($s_t,Q$) policy places a replenishment order of size $Q$ if the inventory level falls below the reorder point $s_t$ and does not place an order otherwise. It is therefore is a special case of ($s_t,Q_t$) in which all $Q_t$'s are equal. We modify the vector space $\mathcal{Q}_t$ introduced in section \ref{sec:sdp-sQt} to explore the ($s_t,Q$) policy.

Let $\dot{\vec{q}}_t\triangleq\langle Q,\ldots,Q\rangle$ be a ($T-t+1$)-dimensional vector of reorder quantities for the ($s_t,Q$) policy and $\dot{\mathcal{Q}}_t$ be a vector space containing all combinations of order quantities $\dot{\vec{q}}_t$. It follows that  $\dot{\mathcal{Q}}_t$ is a subspace of $\mathcal{Q}_t$. For a given $\dot{\vec{q}}_t\in\dot{\mathcal{Q}}_t$, the expected total cost over period $t$ to $T$ when the opening inventory level is $x$ is
\begin{equation}\label{eq:sdp-sQ-Q-It}
V_t(x,\dot{\vec{q}}_t) = \min_{\delta\in\{0,1\}}\mathop{}\{c(\delta Q) + L_t(x+\delta Q)+ \mathbb{E}[V_{t+1}(x+\delta Q - d_t,\dot{\vec{q}}_{t+1})]\},
\end{equation}
and
\begin{equation}\label{eq:sdp-sQ-Q-IT}
V_T(x,\dot{\vec{q}}_T) = \min_{\delta\in\{0,1\}}\mathop{}\{c(\delta Q)+L_T(x+\delta Q)\}
\end{equation}
as the boundary condition. The optimal expected total cost under the ($s_t,Q$) policy with opening inventory level $x$ can be defined as
\begin{equation}
V_0(x) = \min_{\dot{\vec{q}}_1\in\dot{\mathcal{Q}}_1} \mathop{}\{V_1(x,\dot{\vec{q}}_1)\}. \label{eq:sdp-sQ-Q-I0}
\end{equation}
We let the optimal order quantity vector be $\dot{\vec{q}}_t^* \triangleq \langle Q^*,\ldots,Q^* \rangle$. Since $\dot{\mathcal{Q}}_t$ is a subspace of $\mathcal{Q}_t$, Lemma \ref{lemma:Gt-non-increasing} holds for any $\dot{\vec{q}}_t\in\dot{\mathcal{Q}}_t$. The determination of reorder points under ($s_t,Q$) follows the same fashion as ($s_t, Q_t$) policy by Eq.\eqref{eq:sQt-s}. We denote the reorder points associated with $\dot{\vec{q}}_t^*$ as $\dot{\vec{s}}_t^*\triangleq \langle \dot{s}_t^*, \ldots,\dot{s}_T^* \rangle$.

\begin{exmp}[continues=ex1]
Recall the 4--period stochastic lot-sizing problem under Poisson-distributed demand with rates  $d_t = \langle20,40,60,40\rangle$. Under the ($s_t,Q$) policy, the optimal order quantity is $Q^* = 83$ as illustrated by Fig. \ref{fig:ex1.sQ_sdp-Qvscost}. 
The reorder points associated with $\dot{\vec{q}}_1^*$ are determined as $\dot{\vec{s}}_1^* = \langle13, 33, 54, 24\rangle$. Fig. \ref{fig:ex1.sQ-sdp-s1-diff} and \ref{fig:ex1.sQ-sdp-s1-compare} illustrate determining $\dot{s}_1^* = 13$. Note that we apply curves to show the trend of expected costs, while the system is in fact discrete. In Fig. \ref{fig:ex1.sQ-sdp-s1-compare}, a unique sign change of $[\Delta J_1(I_0,\dot{\vec{q}}^*_1) - c(Q^*)]$ is detected between $I_0 = 12$ and 13 and so, by Eq.\eqref{eq:sQt-s}, $I_0 = 13$ is chosen as $\dot{s}^*_1$.
\begin{figure}[!h]
\centering
\includegraphics[width=0.8\textwidth]{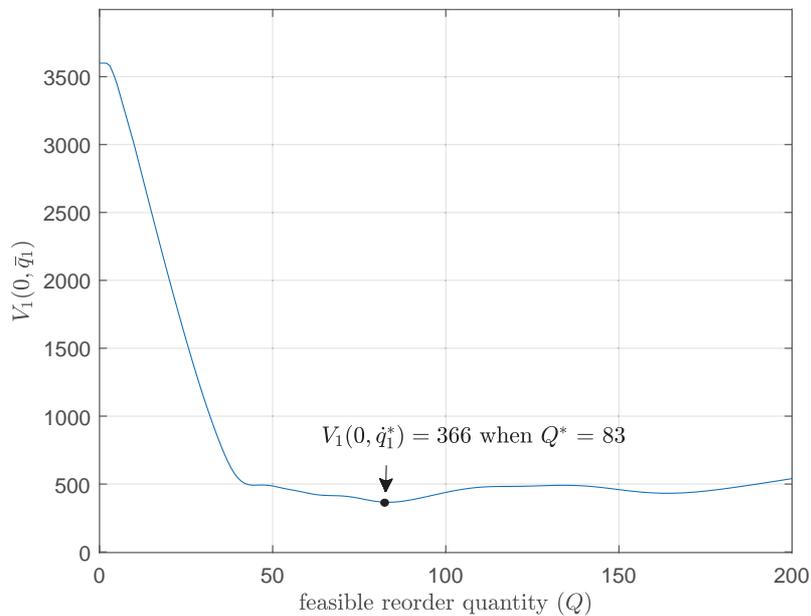}
\caption{$Q^* = 83$ under ($s_t,Q$) policy for Example \ref{ex:classicLarge}.}\label{fig:ex1.sQ_sdp-Qvscost}
\end{figure}
\begin{figure}[H]
    \centering
    \includegraphics[width=0.8\textwidth]{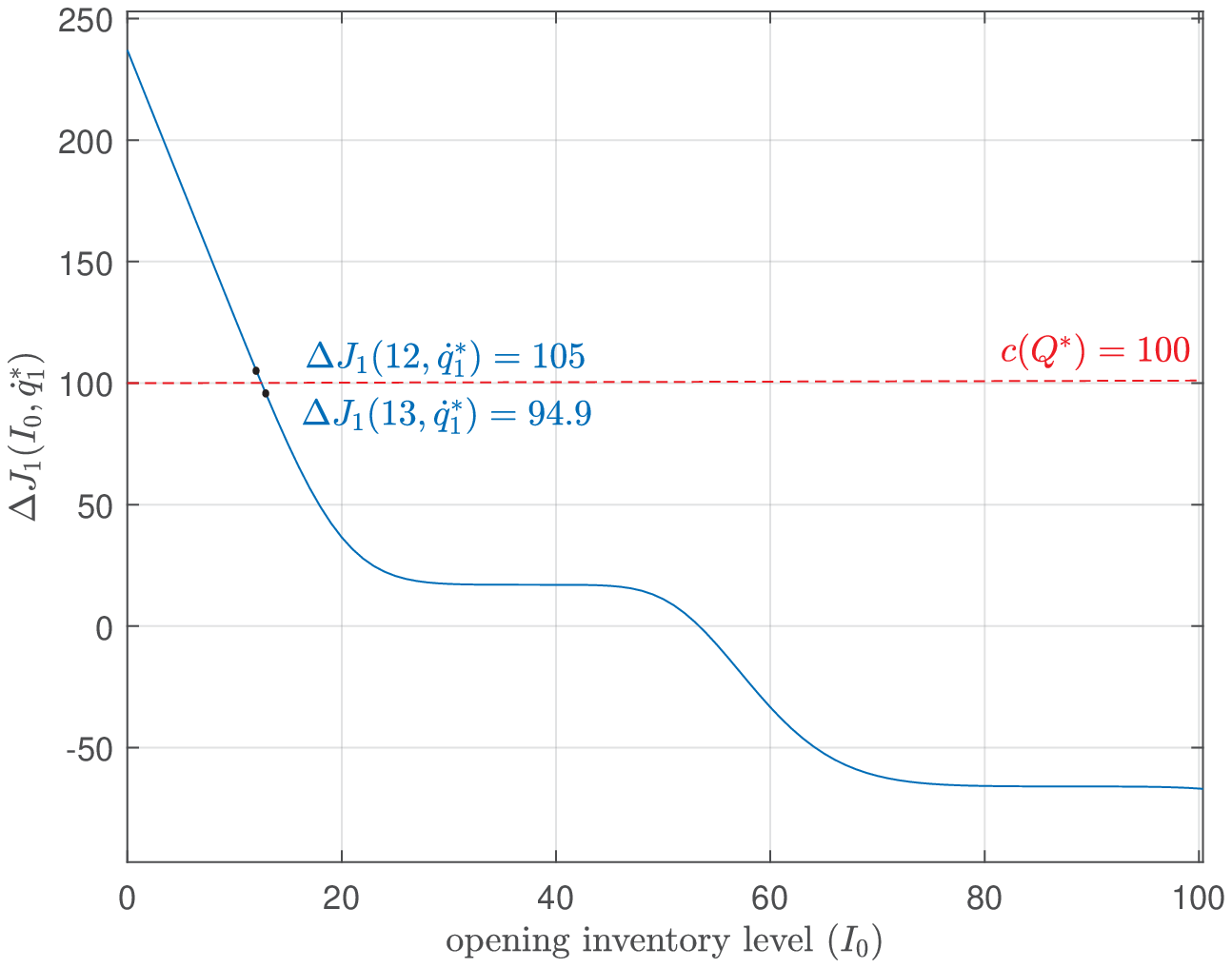}
    \caption{$\dot{s}_1^* = 13$ determined by comparing $\Delta J_1(I_0,\dot{\vec{q}}^*_1)$ and $c(Q^*)$.}\label{fig:ex1.sQ-sdp-s1-diff}
\end{figure}
\begin{figure}[H]
    \centering
    \includegraphics[width=0.8\textwidth]{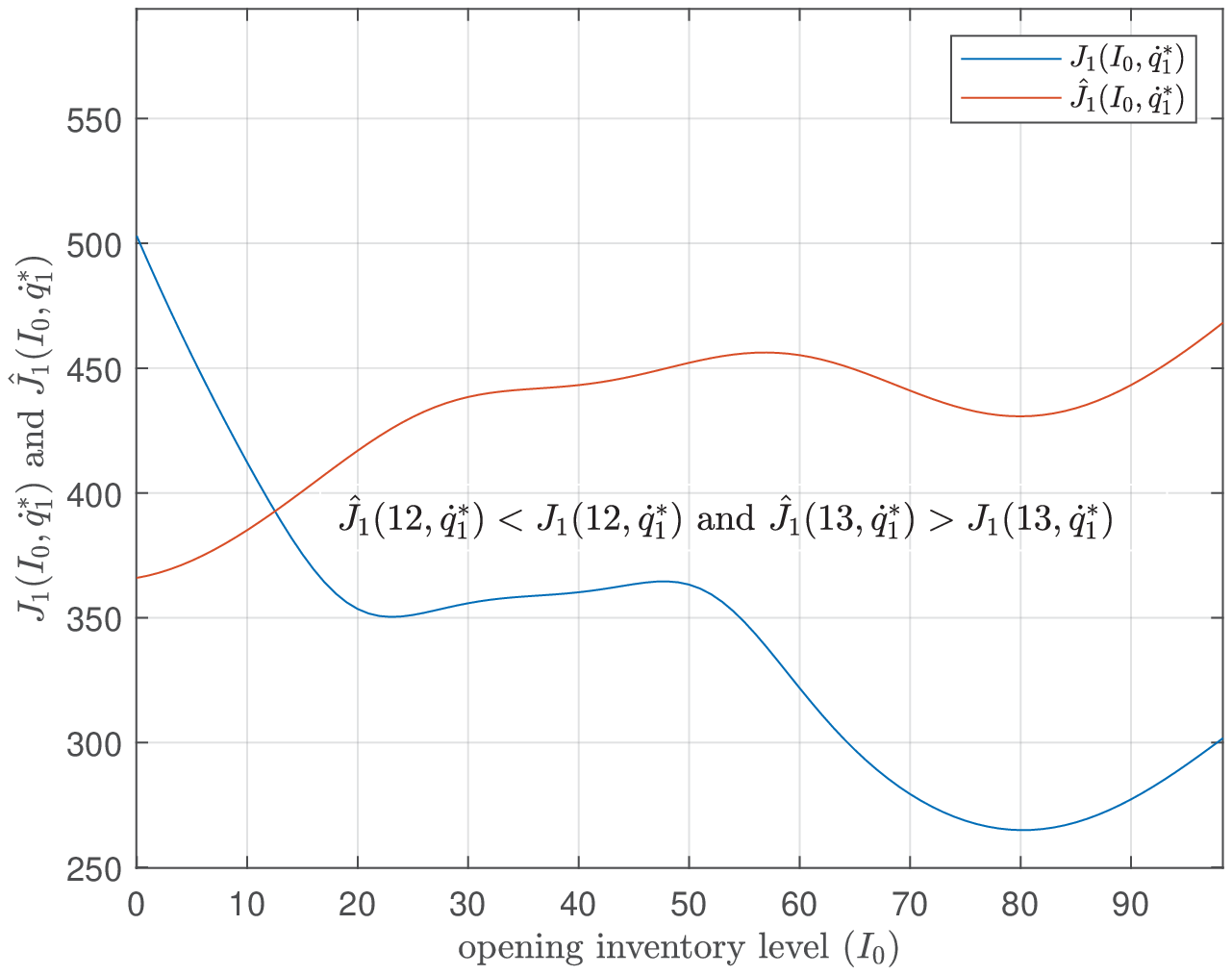}
    \caption{$\dot{s}_1^* = 13$ determined by comparing $J_1(I_0,\dot{\vec{q}}^*_1)$ and $\hat{J}_1(I_0,\dot{\vec{q}}^*_1)$.}
    \label{fig:ex1.sQ-sdp-s1-compare}
\end{figure}
\end{exmp}


\section{A MINLP-based heuristic algorithm for (\textit{$\bm{s_t,Q_t}$})  policy}\label{sec:heuristic}
Optimal ($s_t, Q_t$) and ($s_t, Q$) policies can be obtained by enumerating all possible order quantities and using the stochastic dynamic programming formulations presented in Section \ref{sec:sdps}. However, when the length of planning horizon increases, the enumeration increases exponentially and it becomes impractical to use this method.
In this section, we therefore introduce an effective heuristic to compute near-optimal ($s_t, Q_t$) and ($s_t, Q$) policy parameters in reasonable time. Our heuristic leverages a MINLP approximation of $V_t(\cdot)$ and, similarly to \citet{bookbinder1988strategies}, it comprises two steps: in the first step, we determines a set of near-optimal order quantities; in the second step, we compute the associated reorder points.

\subsection{Step I: Order quantity \textit{$\bm{Q_t}$} of (\textit{$\bm{s_t,Q_t}$}) policy}\label{sec:step1}
We first aim to derive a vector of near-optimal order quantities $\hat{\vec{q}}_t\triangleq\langle\hat{Q}_1.\ldots,\hat{Q}_T\rangle$ for our heuristic ($s_t,Q_t$) policy.
The reader should note that we seek a policy that is near-optimal in terms of expected total cost, not in terms of how close the policy parameters obtained are to the true optimal ones. Therefore, our approximated order quantities and reorder points do not need to be close to the true optimal ones for the $(s_t, Q_t)$ policy, as long as the expected total cost they provide is close enough to the expected total cost of an optimal policy.

Note that if an order is placed in period $t$ under the ($s_t,S_t$) policy, the order quantity is at least $S_{t} - s_{t}$; in fact, if the opening inventory level $I_{t-1} < s_t$ in period $t$, a further $s_t - I_{t-1}$ items will be ordered to ensure the order-up-to level is reached.
In our heuristic ($s_{t},Q_{t}$) policy, we define $\hat{Q}_t\triangleq S_t-s_t$ to be our approximate order quantity in period $t$; and we will denote the vector of approximate order quantities as $\vec{\hat{q}}_t \triangleq \langle \hat{Q}_t, \ldots, \hat{Q}_T \rangle$. While these $\hat{Q}_t$'s may not be optimal, we will compensate for this in Section \ref{sec:step2}, by computing suitable reorder points that are tailored for these approximate order quantities.

Of course, to compute $\hat{Q}_t$, we need optimal or near-optimal values of parameters $s_t$ and $S_t$ of the ($s_t,S_t$) policy. To compute near-optimal $s_t$ and $S_t$ values for large-scale problems, we leverage the approach introduced by \citet{xiang2018computing}. For the sake of completeness, the model we adopted is presented in Appendix \ref{sec:App.MINLP-Xiang}. 

\subsection{Step II: Reorder point \textit{$\bm{s_t}$} of (\textit{$\bm{s_t,Q_t}$}) policy}\label{sec:step2}
Since approximate order quantities $\hat{\vec{q}}_t$ are a lower bound for order quantities observed under an ($s_t,S_t$) policy, we cannot directly use the reorder points from the optimal ($s_t,S_t$) policy as the reorder points for a heuristic ($s_t,Q_t$) policy. To compensate for the under-estimation in the order quantities, we need higher reorder points.

For a given vector $\hat{\vec{q}}_t$ of approximate order quantities, we may compute the associated optimal reorder points by using an SDP formulation.
This would be relatively straightforward for Poisson demand, but would require a discretisation step for continuous demand distributions.
In order to provide a framework that can be applied to Poisson, normal, and possibly other continuous demand distributions, we modify the model in \citep{xiang2018computing} to capture the characteristics of an ($s_t,Q_t$) and provide an approximation $\mathcal{J}_t(x,\hat{\vec{q}}_t)$ of $J_t(x,\hat{\vec{q}}_t)$ that can be used in Eq.\eqref{eq:sQt-s} to compute near-optimal reorder points. This model is named `Model \ref{eq:modelII}'.
Let $\mathcal{J}_t(x,\hat{\vec{q}}_t)$ be our approximation of $J_t(x,\hat{\vec{q}}_t)$ for the set of near-optimal order quantities $\hat{\vec{q}}_t$ computed in Section \ref{sec:step1}.
{\small{
\begin{alignat}{2}
\mathcal{J}_t(x,\hat{\vec{q}}_t) =\min\quad & h\tilde{H}_t + b \tilde{B}_t + \sum_{k=t+1}^T[h \tilde{H}_k + b \tilde{B}_k + c(\delta_k \hat{Q}_k)], &\label{eq:MINLP-obj-step2}\\
\mbox{s.t.}\quad
&\delta_t = 0,&\quad& \label{eq:J-delta=0}\\
&\tilde{I}_t + \tilde{d}_t = \tilde{I}_{t-1},&\quad& \label{eq:J-balance-t} \\
&\delta_k=0 \rightarrow \tilde{I}_k + \tilde{d}_k - \tilde{I}_{k-1} = 0, &\quad& k = t+1,\ldots, T,\label{eq:J-balance-no}\\
&\delta_k=1 \rightarrow \tilde{I}_k + \tilde{d}_k - \tilde{I}_{k-1} = \hat{Q}_k , &\quad& k = t+1,\ldots, T,\label{eq:J-balance-order}\\
& \sum\nolimits_{j=t}^k P_{jk}=1, &\quad& k=t+1,\ldots,T, \label{eq:J-sumReview}\\
& P_{jk} \geq \delta_j - \sum\limits_{r=j+1}^k \delta_r, & \quad & k=t,\cdots,T\text{ and } j = t,\ldots, k,\label{eq:J-latestReview}\\
& P_{jk}=1 \rightarrow \tilde{H}_k = \hat{\mathcal{L}}(\tilde{I}_k+\tilde{d}_{jk},d_{jk}), &\quad& k=t,\ldots,T \text{ and }j = t,\ldots, k,\label{eq:J-H}\\
& P_{jk}=1 \rightarrow \tilde{B}_k = \mathcal{L}(\tilde{I}_k+\tilde{d}_{jk},d_{jk}), &\quad& k=t,\ldots,T \text{ and }j = t,\ldots, k,\label{eq:J-B}\\
& \tilde{H}_k, \tilde{B}_k \geq 0, P_{jk}, \delta_k \in\{0,1\},&\quad& k = t,\ldots, T \text{ and }j = t,\ldots, k.
\end{alignat}\label{eq:modelII}}}
Let $\tilde{H}_k$ and $\tilde{B}_k$ denote the expected positive inventory and back-ordered levels at the end of period $k$, respectively; their values are computed by following the piecewise-linear approximation strategy in \citet{rossi2015piecewise}, which is based on the first-order loss function $\mathcal{L}$ and its complement $\hat{\mathcal{L}}$. We discuss in detail on the loss function and piecewise-linear approximation under non-stationary demand of Poisson distribution in Appendix \ref{sec:App.PoissonPiecewise}.

In line with \citep{tarim2006modelling}, $\delta_k$ is a binary variable that takes value 1 if and only if an order is placed in period $k$, while $P_{jk}$ is a binary variable that takes the value 1 if and only if the most recent inventory review\footnote{An inventory review is a point in time at which we observe the inventory level, which therefore becomes a known quantity.} before period $k$ took place at the beginning of period $j$; note that variable $P_{jk}$ allows us to properly account for demand variance while computing the first-order loss function.
In the model above, the objective function $\mathcal{J}_t(x,\hat{\vec{q}}_t)$ approximates the expected total cost over horizon ($t,T$) with no order in period $t$. In contrast to \citeauthor{xiang2018computing}'s model, the order quantities in this revised model are no longer decision variables, but a set of near-optimal policy parameters $\hat{\vec{q}}_t$ obtained in Section \ref{sec:step1}. 
We add constraint \eqref{eq:J-delta=0} to ensure that no order is placed in the first period of the planning horizon ($t,T$). We also modify the flow balance in period $t$ as constraint \eqref{eq:J-balance-t}, where $\tilde{I}_t$ denotes the expected closing inventory of period $t$. The other constraints remain as in \citep{xiang2018computing}.


Since $J_t(x,\vec{q}_t)$ is approximated as $\mathcal{J}_t(x,\hat{\vec{q}}_t)$, the near-optimal reorder point $\hat{s}_t$ can be determined, following Eq.\eqref{eq:sQt-s}, as
\begin{equation}
\hat{s}_t = \inf\{x:\mathop{} \Delta \mathcal{J}_t(x,\hat{\vec{q}}_t) < K+z\hat{Q}_t\},\label{eq:stQt-MINLP-s-hard}
\end{equation}
or as the minimum value of $x$ such that
\begin{equation}\label{eq:stQt-MINLP-s-hard-int}
\Delta \mathcal{J}_t(x,\hat{\vec{q}}_t) < K+z\hat{Q}_t
\end{equation}
for discrete inventory levels, where $\Delta \mathcal{J}_t(x,\hat{\vec{q}}_t) \triangleq \mathcal{J}_t(x,\hat{\vec{q}}_t) - \mathcal{J}_t(x+\hat{Q}_t,\hat{\vec{q}}_t)$. Note that, there is no guarantee of monotonicity for $\Delta \mathcal{J}_t(x,\hat{\vec{q}}_t)$ in $x$ since the piecewise linearisation produces errors; our model applies the optimal partitioning strategy to maintain a minimum error \citet[Thm.~11]{rossi2014piecewise}. We denote the vector of near-optimal reorder points associated with $\vec{\hat{q}}_t$ as $\vec{\hat{s}}_t\triangleq\langle \hat{s}_t,\ldots, \hat{s}_T\rangle$.

\subsection{A binary search approach to approximate the reorder points $\bm{s_t}$}\label{sec:algorithm}
A line search for $\hat{s}_t$ following Eq.\eqref{eq:stQt-MINLP-s-hard} may be too time-consuming for large-scale instances. This subsection introduces a heuristic algorithm to approximate $\hat{s}_t$ and reduce computational complexity.

The algorithm applies a binary search on $\Delta \mathcal{J}_t(x,\vec{\hat{q}}_t)$ with $\vec{\hat{q}}_t$ known as an input. For any period $t$, input opening inventory level $x_0$ and given step-size $w$ ($w>0$) define an interval of inventory level $[x_0,\mathop{} x_0 + w]$, which maps to $[\Delta \mathcal{J}_t(x_0+w ,\hat{\vec{q}}_t),\mathop{} \Delta \mathcal{J}_t(x_0,\hat{\vec{q}}_t)]$. The binary search halves the length of the interval in each iteration until $\hat{s}_t$ is detected according to Eq.\eqref{eq:stQt-MINLP-s-hard}. If the initial interval does not span the point at which the sign of $\Delta \mathcal{J}_t(x,\hat{\vec{q}}_t) - K - z\hat{Q}_t$ changes, we renew $[x_0,\mathop{} x_0 + w]$ by panning it $w$ units to the left if $\Delta \mathcal{J}_t(x_0+w,\hat{\vec{q}}_t) <  K + z\hat{Q}_t$ or to the right, otherwise; and then proceed with the binary search.

We present the following algorithm for integer inventory levels. One can extend it to discrete systems with any interval between two adjacent inventory levels. For integer inventory levels, the algorithm terminates if a pair of inventory levels $x$ and $x+1$ are found such that $\Delta\mathcal{J}_t(x,\hat{\vec{q}}_t)\leq K+z\hat{Q}_t\leq\Delta\mathcal{J}_t(x+1,\hat{\vec{q}}_t)$, and then $\hat{s}_t=x+1$. The procedure in detail is as follows.

\vspace{-0.1cm}
\alglanguage{pseudocode}
\begin{algorithm}[H]\renewcommand\baselinestretch{1.2}\selectfont
\small
\caption{Computing the reorder points $\hat{s}_t$ associated with $\vec{\hat{q}}_t$.}
\label{Algorithm}
\begin{algorithmic}[1]
\State \textbf{Input:} demand rates $\tilde{d}_t$;
cost parameters ($K$, $z$, $h$, $b$);
the step-size $w$;
an opening inventory $x_0$;
order quantities $\vec{\hat{q}}_t$.
\State \textbf{Output:} reorder point $\hat{s}_t$ associated with $\vec{\hat{q}}_t$.
\For {$t = 1 \to T$}
    \State Compute the ordering cost of placing an order $\mathcal{J}_0 = K + z\hat{Q}_t$;
    \State $x_l = x_0$ and $x_r = x_0 + w$;
    \State compute $\mathcal{J}_l = \Delta \mathcal{J}_t(x_l,\hat{\vec{q}}_t)$ and $\mathcal{J}_r = \Delta \mathcal{J}_t(x_r,\hat{\vec{q}}_t)$ with $\hat{Q}_t$;
     \If {$\mathcal{J}_l>\mathcal{J}_0>\mathcal{J}_r$}
        \State $x_m = \lfloor\frac{x_l+x_r}{2}\rfloor$ and $\mathcal{J}_m = \Delta \mathcal{J}_t(x_m,\hat{\vec{q}}_t)$;
        \If {$\mathcal{J}_m > \mathcal{J}_0$}
            \If{$\Delta \mathcal{J}_t(x_m+1,\hat{\vec{q}}_t)<\mathcal{J}_0$}
                \State output $\hat{s}_t = x_m$;
                \Else \hspace{0.5em}$x_l = x_m$, $x_r = x_r$, and repeat lines 6 -- 20;
            \EndIf
            \Else
            \If {$\Delta \mathcal{J}_t(x_m-1,\hat{\vec{q}}_t)>\mathcal{J}_0$}
                \State output $\hat{s}_t = x_m-1$;
                \Else \hspace{0.5em} $x_l = x_l$, $x_r = x_m$, and repeat lines 6 -- 20;
            \EndIf
        \EndIf
    \EndIf
\EndFor
\end{algorithmic}
\end{algorithm}

\begin{exmp}[continues=ex2]
Recall the 4--period stochastic lot-sizing problem under Poisson-distributed demand with rates $d_t = \langle2, 1, 5, 3\rangle$. Applying 20 partitions in the piecewise linearisation approximation, $\vec{\hat{q}}_1 = \langle3, 4, 9, 5 \rangle$ approximates $J_1(I_0,\vec{q}_1^*)$ as shown in Fig. \ref{fig:ex1.sQt-MINLP-sdp-compare} for $I_0\in[-4,14]$. The curves are plotted to demonstrate the difference between $J_1(I_0,\vec{q}_1^*)$ and $\mathcal{J}_t(I_0)$, while the system is in fact discrete.
\begin{figure}[h!]
\centering
\includegraphics[width=0.8\textwidth]{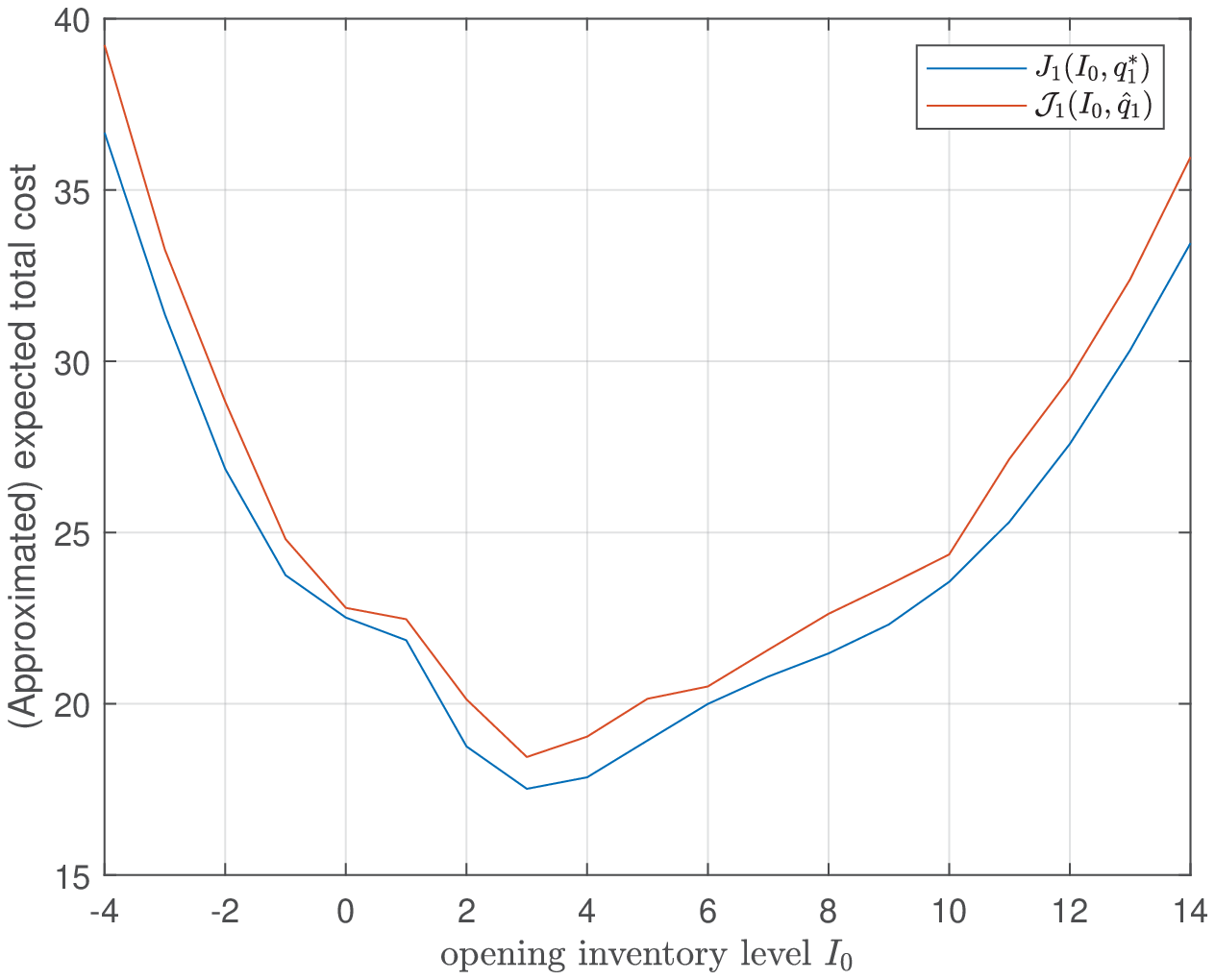}
\caption{Plot of $J_1(I_0,\vec{q}_1^*)$ and $\mathcal{J}_1(I_0,\hat{\vec{q}}_1)$.} \label{fig:ex1.sQt-MINLP-sdp-compare}
\end{figure}
\noindent Table \ref{tab:sQt-MINLP-exmp-compareParameter} compares the ($s_t,Q_t$) policy parameters attained by the SDP in Section \ref{sec:sdp-sQt} and the heuristic for a zero initial inventory level. 
\vspace{-0.4cm}
\begin{table}[H]\renewcommand\arraystretch{1}\selectfont\small
\caption{Policy parameters of Example \ref{ex:smallPoisson_stQt} computed by SDP and heuristic under ($s_t,Q_t$) policy.}\label{tab:sQt-MINLP-exmp-compareParameter}
\centering
\begin{tabular}{rp{0.5cm}<{\centering} p{0.5cm}<{\centering} p{0.5cm}<{\centering} p{0.5cm}<{\centering}  p{0.5cm}<{\centering} p{0.5cm}<{\centering} p{0.5cm}<{\centering} p{0.5cm}<{\centering}}
\hline
\multirow{2}{*}{} & \multicolumn{4}{c}{$\hat{Q}_t$} & \multicolumn{4}{c}{$\hat{s}_t$}\\
\cmidrule(r){2-5}\cmidrule(r){6-9}
$t$& 1 & 2 & 3 & 4 & 1 & 2 & 3 & 4\\
\hline
{\bf SDP}&  3  & 3     & 8     & 5& 1     &0      & 4     & 1\\
{\bf Heuristic}&  3   &4      &9      &5  &  1     &-2      & 4     & 0\\
\hline
\end{tabular}
\end{table}
\noindent Taking $G_1(0) = 21.8$ as a benchmark, the optimality gaps of the ($s_t, Q_t$) determined policy by SDP and our heuristic, relative to the ($s_t,S_t$) policy are showed in Table \ref{tab:sQt-MINLP-exmp-compareGap}. We note that the ($s_t,Q_t$) policy produces large optimality gaps  in Example \ref{ex:smallPoisson_stQt}, where ETC values are small, while the approximation accuracy of the heuristic $(23.1 - 22.5)/22.5\times 100\% = 2.67\%$ is acceptable. We will extend the computation in Section \ref{sec:computations} to investigate ($s_t,Q_t$) policy performs on extensive instances.
\vspace{-0.4cm}
\begin{table}[htb]\renewcommand\arraystretch{1}
    \caption{Expected total cost (ETC) and optimality gap (OG) of Example \ref{ex:smallPoisson_stQt} by SDP and heuristic under the ($s_t,Q_t$) policy.}
    \label{tab:sQt-MINLP-exmp-compareGap}
    \centering
    \begin{tabular}{rcc}
    \hline
        & ETC & OG($\%$)\\
        \hline
        {\bf SDP} & 22.5& $3.33$\\
        {\bf Heuristic}&23.1&$5.93$\\
        \hline
    \end{tabular}
\end{table}
\end{exmp}

\subsection{Approximation of (\textit{$\bm{s_t}$,$\bm{Q}$}) policy parameters}\label{sec:MINLP-stQ-difference}
For the ($s_t,Q$) policy, a direct way to approximate the order quantity is to simplify model in Appendix \ref{sec:App.MINLP-Xiang} by replacing  $Q_t$ with $Q$ and then follow the steps in Sections \ref{sec:step1} and \ref{sec:step2}; however, this is found to produce large optimality gaps in terms of the expected total cost.

Following the line of reasoning illustrated in Section \ref{sec:step1} for ($s_t,Q_t$), one can derive a single order quantity in period 1 as $S_1 - I_0$ for a known opening inventory $I_0$.
However, a high value for $I_0$ may result in a low order quantity imposed over a long period.
In our heuristic ($s_t,Q$) policy, we define $\hat{Q} \triangleq S_1$ to be our approximate order quantity for horizon ($1,T$); and we denote the vector of approximate order quantities as $\hat{\vec{q}} \triangleq \langle \hat{Q},\ldots,\hat{Q} \rangle$.
The reorder points are adjusted to compensate for the over-estimation for cases with high opening inventory levels.

The determination of reorder points $\hat{\vec{s}}_t$ associated with order quantity $\hat{Q}$ follows the same procedure proposed in Section \ref{sec:step2} for ($s_t,Q_t$). We apply Model \ref{sec:step2} with $\hat{Q}$ to obtain the approximated expected cost over horizon ($t,T$) when no order is placed in $t$, denoted as $\mathcal{J}_t(x,\hat{\vec{q}})$, and we apply our previously introduced heuristic algorithm on the function $\Delta \mathcal{J}_t(x,\hat{\vec{q}})$ to determine $\hat{s}_t$.

\begin{exmp}[continues=ex1]
Applying 20 partitions in the piecewise linearisation approximation, Fig. \ref{fig:ex1.sQ-Q-MINLP-sdp-compare} approximates $J_1(I_0,\dot{\vec{q}}^*_1)$ by $\mathcal{J}_1(I_0,\hat{\vec{q}})$. Similarly, the inventory system is discrete, while we apply curves to demonstrate the difference.
\begin{figure}[htb]
\centering
\includegraphics[width=0.8\textwidth]{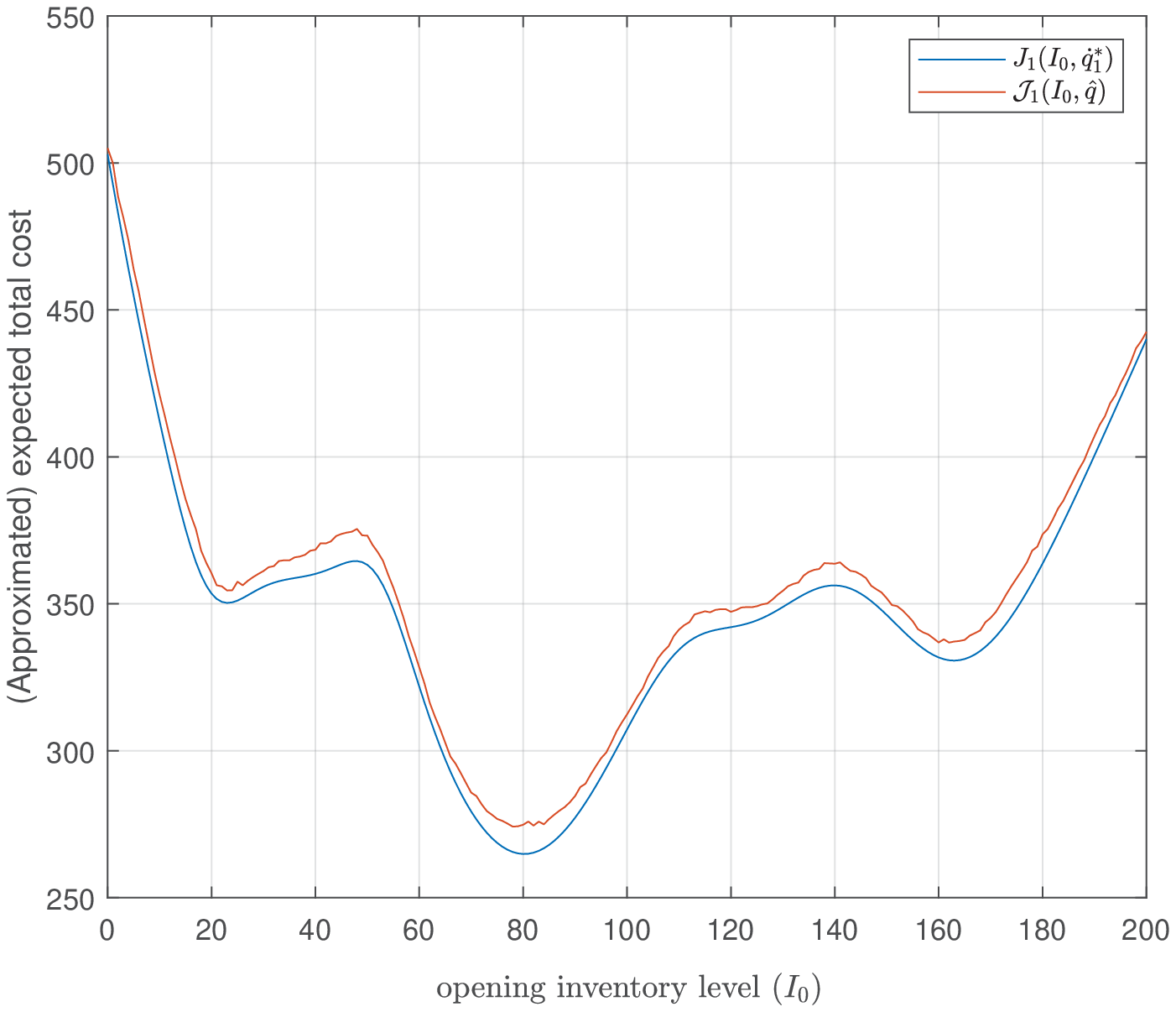}
\caption{Plot of $J_1(I_0,\dot{\vec{q}}_1^*)$ and $\mathcal{J}_1(I_0,\hat{\vec{q}})$.}\label{fig:ex1.sQ-Q-MINLP-sdp-compare}
\end{figure}
\noindent For a zero initial inventory level, Table \ref{tab:sQ-Q-MINLP-exmp-compare} compares the policy parameters computed by the SDP and approximation under the ($s_t,Q$) policy.
\begin{table}[htb]\renewcommand\arraystretch{1}\selectfont\small
\caption{Policy parameters of Example \ref{ex:classicLarge} computed by SDP and heuristic under the ($s_t,Q$) policy.}\label{tab:sQ-Q-MINLP-exmp-compare}
\centering
\begin{tabular}{rp{0.2cm}<{\centering} p{0.00001cm}<{\centering} p{0.5cm}<{\centering}  p{0.5cm}<{\centering} p{0.5cm}<{\centering} p{0.5cm}<{\centering}}
\hline
\multirow{2}{*}{} & \multicolumn{2}{c}{$\hat{Q}$} & \multicolumn{4}{c}{$\hat{s}_t$}\\
\cmidrule(r){2-3}\cmidrule(r){4-7}
$t$ &-- & & 1 & 2 & 3 & 4\\
\hline
{\bf SDP}&  83& & 13 & 33 & 54 & 24\\
{\bf Heuristic}&  84& & 14 &34 & 55 & 24\\
\hline
\end{tabular}
\end{table}
\noindent Taking $G_1(0) = 481$ as the benchmark, Table \ref{tab:sQt-MINLP-exmp-compareGap} summarises the optimality gaps of the ($s_t, Q$) policy by SDP and the heuristic. The approximation accuracy $(504 - 502 )/502\times 100\% = 0.398\%$ behaves well. We discuss the performance of the ($s_t,Q$) policy in detail in the next section.
\begin{table}[htb]\renewcommand\arraystretch{1}\selectfont\small
    \caption{Expected total cost (ETC) and optimality gap (OG) of Example \ref{ex:classicLarge} by ($s_t,Q$) with SDP and heuristic.}
    \label{tab:sQt-MINLP-exmp-compareGap}
    \centering
    \begin{tabular}{rcc}
    \hline
        & ETC & OG($\%$) \\
        \hline
        {\bf SDP} & 503& $4.57$\\
        {\bf Heuristic}&505&$4.99$\\
        \hline
    \end{tabular}
\end{table}
\end{exmp}


\section{Computational analysis}\label{sec:computations}
This section presents a computational analysis to evaluate ($s,Q$)-type policies under non-stationary stochastic demand. The analysis considers both the stochastic dynamic programming formulations and our heuristics for the ($s_t, Q_t$) and ($s_t, Q$) policies.
In Section \ref{sec:6}, we consider a test set comprising small 6-period instances; we investigate the performances of optimal ($s, Q$)-type policies against optimal non-stationary ($s,S$) policy, and we evaluate the difference between optimal ($s,Q$) and heuristic ($s,Q$) policies.
In Section \ref{sec:25}, we consider a large test set comprising 25-period instances; we investigate the performance of ($s,Q$)-type heuristics versus the optimal non-stationary ($s,S$) policy; we also compare the performance between our ($s,Q$)-type heuristics and another existing static-dynamic uncertainty heuristic, namely the ($R_t,S_t$) policy discussed in \citep{rossi2015piecewise}.

We name the optimal policy for the stochastic lot-sizing problem, which takes an ($s,S$) form, ($s_t,S_t$)-SDP.
In our experiment we consider two variants of the ($s,Q$) policy: the ($s_t,Q_t$) policy, and the ($s_t,Q$) policy; presented in Section \ref{sec:sdp-sQt} and Section \ref{sec:sdp-sQ-Q}, respectively.
For each variant, we discuss results for the optimal SDP formulation, named ($s_t,Q_t$)-SDP and ($s_t,Q$)-SDP, respectively; and results for our MINLP heuristics formulations presented in Section \ref{sec:heuristic}, named ($s_t,Q_t$)-Heuristic and ($s_t,Q$)-Heuristic, respectively.
We apply 10 partitions in the piecewise approximation for both heuristics.
We simulate each test instance with the policy parameters obtained from the heuristics and derive the average total cost of 500,000 simulation runs.

For each approach, we always use the optimal ($s,S$) policy as a benchmark. Approaches are compared in terms of their expected total cost (ETC) percent optimality gap computed as $100\times(\text{ETC}_2 - \text{ETC}_1)/\text{ETC}_1$, where $\text{ETC}_1$ is the expected total cost of the optimal non-stationary ($s,S$) policy, and $\text{ETC}_2$ is the expected total cost of the other approach benchmarked.
We set a zero initial inventory for all test instances and test the robustness of heuristics for ($s,Q$)-type policies.

In our numerical study, we consider ten expected demand patterns: two life cycle patterns, one moves from the launch stage to maturity via a growth (LCY1) and the other moves from the growth stage through maturity and into decline (LCY2); two sinusoidal patterns, one with stronger (SIN1) and the other with weaker (SIN2) oscillations; a stationary demand pattern (STAT); a random demand pattern (RAND); and lastly, 4 empirical patterns derived according to \citep{strijbosch2011interaction}.

All computations are performed by a 4.0 (1.90+2.11) gigahertz Intel(R) Core(TM) i7-8650U CPU with 16.0 gigabytes of RAM in JAVA 1.8.0\_201.

\subsection{A test set with 6-period Poisson-distributed demand}\label{sec:6}

The first test set involves 60 instances over a 6-period planning horizon in which the demand follows a non-stationary Poisson distribution. Our aim is twofold: first, we aim to investigate the performances of optimal ($s, Q$)-type policies obtained via SDP against the optimal non-stationary ($s,S$) policy; second we aim to evaluate the difference between optimal ($s,Q$) and heuristic ($s,Q$) policies.

\begin{figure}[h!]
\centering
\includegraphics[width=0.9\textwidth]{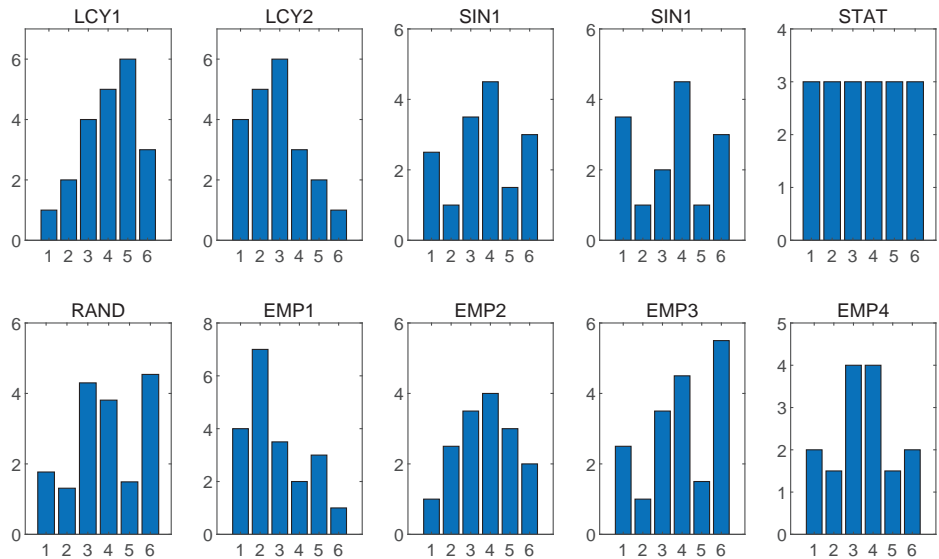}
\caption{Demand patterns of 6-period instances.}\label{fig:6}
\end{figure}

We assume the maximum order quantity is 9, which allows us to enumerate all combination of order quantities for the ($s_t,Q_t$) policy by stochastic dynamic programming. The problems in this test set are designed with very small mean demands $\lambda_t$, as illustrated in Fig. \ref{fig:6}.
The values of $\lambda_t$ are set to be between 1 and 7 in all cases which allows variation in the optimal values of $Q_t$ and ensures that the optimal order quantity is never as high as 9 in any period. The problem coefficients are considered over $z \in \{0,1\}$ and the three sets of $K$ and $b$ shown in Table \ref{tab:6parameters} with different ratios of $K$ to $b$. Holding cost is set as $h = 1$ for all instances.

\begin{table}[H]\renewcommand\arraystretch{0.8}\selectfont\small
\centering
\caption{parameter groups of fixed ordering cost ($K$) and penalty cost ($p$)}\label{tab:6parameters}
\begin{tabular}{cccc}
\toprule
set &	$K$ & $b$ & ratio\\
\midrule
1 &	5 &	3 & 1.67\\
2 &	10 &	3 & 2.00\\
3 &	10 &	7 & 1.43\\
\bottomrule
\end{tabular}
\end{table}

Table \ref{tab:6gap} reports, for each approach considered, the optimality gaps observed against the optimal ($s_t,S_t$) policy. The results for ($s_t,Q_t$)-SDP and ($s_t,Q$)-SDP give the exact optimality gaps for these policies against optimal ($s_t,S_t$) policy, which are on average $1.91\%$ and $3.61\%$ respectively. In detail, ($s_t,Q_t$)-SDP performs better than ($s_t,Q$)-SDP in every individual demand pattern; and ($s_t,Q$)-SDP is dominated by ($s_t, Q_t$)-SDP even in the case of a stationary demand pattern. In view of cost parameters, there is no obvious relation between optimality gaps and the variation in demand patterns or in the ratio of $K$ to $b$. Optimality gaps also remain consistent when the unit cost is changed. On the other hand, the increase in penalty cost results in a small increase in the optimality gap for both ($s_t,Q_t$)-SDP ($1.43\%$ to $1.49\%$) and ($s_t,Q$)-SDP ($2.92\%$ to $3.22\%$).

For ($s_t,Q_t$)-Heuristic and ($s_t,Q$)-Heuristic, we found average differences of $0.85\%$ and $1.05\%$. The largest average difference arises under demand pattern EMP3 ($1.04\%$) for ($s_t,Q_t$)-Heuristic and RAND ($1.75\%$) for ($s_t,Q$)-Heuristic. We conclude that the difference between SDP and the heuristic approach is generally low.

\begin{table}[H]\renewcommand\arraystretch{0.8}\selectfont\small
\centering
\caption{Average percent ETC optimality gap over our 6-period test set under different demand patterns and pivoting parameters.
}\label{tab:6gap}
\begin{tabular}{cp{1.5cm}<{\centering}p{1.5cm}<{\centering}p{1.5cm}<{\centering}p{1.5cm}<{\centering}}
\toprule
Problem Settings & ($s_t,Q_t$)-SDP   & ($s_t,Q_t$)-Heuristic   &   ($s_t,Q$)-SDP   & ($s_t,Q$)-Heuristic\\
\midrule
\textbf{demand pattern} &&&\\
LCY1&1.96&2.60&2.55&3.30\\
LCY2&2.70&3.60&5.37&6.11\\
SIN1&1.95&2.89&3.96&4.80\\
SIN2&2.13&3.04&3.18&4.75\\
STAT&1.54&2.41&2.45&4.00\\
RAND&1.17&2.02&3.12&4.86\\
EMP1&1.98&2.87&3.98&5.33\\
EMP2&2.32&2.94&3.56&4.44\\
EMP3&1.13&2.17&3.11&3.66\\
EMP4&2.21&3.11&4.80&5.39\\
\textbf{unit cost}&&&\\
0&2.03&2.93&3.83&5.15\\
1&1.79&2.59&3.38&4.18\\
\textbf{set}&&&\\
1&2.81&3.76&4.67&5.76\\
2&1.43&2.29&2.92&3.96\\
3&1.49&2.23&3.22&4.28\\
\midrule
\textbf{Average}&1.91&2.76&3.61&4.66\\
\bottomrule
\end{tabular}
\end{table}

\subsection{A test set with 25-period Normally-distributed demand}\label{sec:25}
We extend the planning horizon to 25 periods.
The purpose of implementing this test set is twofold. First we aim to investigate the performance of ($s,Q$)-type heuristics versus the optimal non-stationary ($s,S$) policy for larger instances; second, we aim to compare the performance between ($s,Q$)-type heuristics and the non-stationary ($R,S$) policy introduced in \citep{rossi2015piecewise}, which we name ($R_t,S_t$)-Heuristic.

Since the computation of piecewise linearisation parameters consumes a large amount of computation time for large non-stationary demand following a Poisson distribution, in what follows we will focus on normally distributed demand patterns, for which \citet{rossi2014piecewise} present precomputed optimal partitioning coefficients.

\begin{figure}[h!]
\centering
\includegraphics[width=0.9\textwidth]{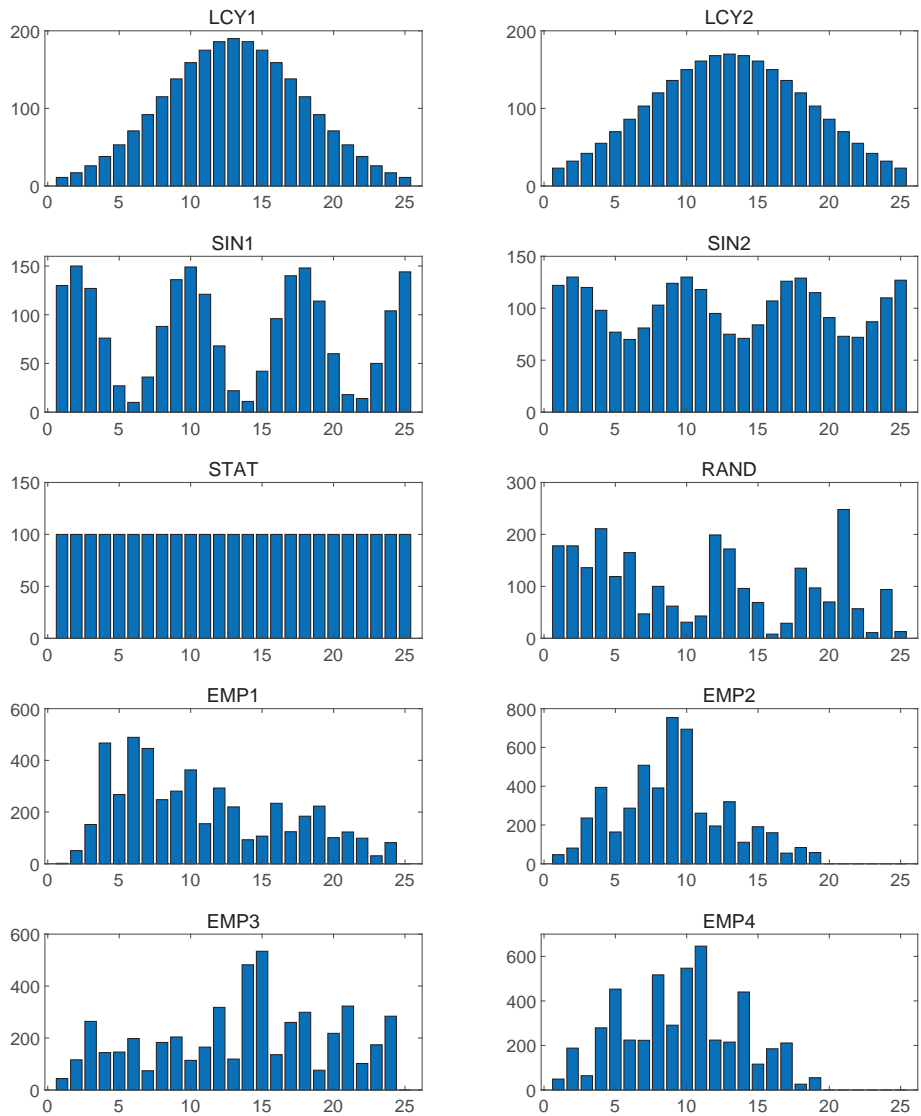}
\caption{Demand patterns of 25-period instances.}\label{fig:25}
\end{figure}

We refer to the 25-period instances in \citep{xiang2018computing}. The demand $d_t$ in each period $t$ is assumed to be a normally distributed random variable with known mean $\tilde{d}_t$ and standard deviation $\sigma_t = \rho\cdot \tilde{d}_t$, where $\rho$ denotes the coefficient of variation of the demand, which remains fixed over time as prescribed in \citep{bollapragada1999simple}; demands are assumed to be independent among each other. We allow the standard deviation parameter $\rho$ to vary over $\rho \in \{0.1, 0.2, 0.3\}$. Demand patterns are illustrated in Figure \ref{fig:25}. Other problem parameters are $K\in\{500, 1000, 1500\}$; $b \in \{5, 10, 20\}$; $z \in \{0,1\}$; and $h = 1$.

The reader should note that, since stochastic dynamic programming is pseudo-polynomial, an increase in the average value of the demand or of its standard deviation will lead to a dramatic increase in the state space and hence of computational times \citep{dural2019benefit}. The ($s_t,Q$)-SDP can be implemented by bounding the inventory level, while it is no longer possible to compute ($s_t,Q_t$)-SDP within a reasonable time for normal demand or large planning horizons such as 25.

Table \ref{tab:25gaps} reports average optimality gaps for our 25-period instances. For the ($s_t,Q_t$)-Heuristic, the average optimality gap of ETC is $2.31\%$, which is similar to the result obtained for the 6-period test problems. The optimality gap exhibits similar trends with the penalty cost and the unit cost, while the gap increases with penalty cost and decreases when the unit cost is increased. For the normal distribution, the increase of the standard deviation parameter $\rho$ reduces the optimality gap, which suggests the ($s_t,Q_t$) policy performs slightly better when the demand standard deviation is higher.

For ($s_t,Q$)-Heuristic, once more, as with 6-period test set, we observe that the ($s_t,Q$)-Heuristic is not satisfactory. The average optimality gap now increases up to $11.5\%$; and for an individual demand pattern, the optimality gap reaches $25.9\%$. We also cross-validated results against optimal ($s_t,Q$) parameters obtained via SDP, to ensure the accuracy of the result, but found that the optimality gap remained as large as $10.5\%$ on average. This confirms that not just the approximation, but the policy itself performs poorly. We believe that, when the length of planning horizon increases, under non-stationary demand the single order quantity $Q$ in ($s_t,Q$) policy cannot properly hedge against demand, and thus it produces substantially higher expected cost than other policies that provide more flexibility. It should be noted that the maximum optimality gaps observed for ($s_t,Q$)-SDP ($24.9\%$ and $23.8\%$) concern empirical demand patterns with a series of 0 demand. A single order quantity for all periods causes either a large amount of holding cost for 0-demand periods or penalty cost for large-demand periods. Despite the unsatisfactory performance of the ($s_t,Q$) policy, it is worth noting that the results show the same trends with respect to $\rho$, $b$ and $z$ as the ($s_t,Q_t$) policy.

The optimality gaps of ETC by ($R_t,S_t$)-Heuristic is found as $2.90\%$, which are larger than those observed for ($s_t,Q_t$)-Heuristic over all demand patterns and pivoting parameters. As a result, we conclude that in the context of our test set the ($s_t,Q_t$) is better than ($R_t,S_t$) policy in terms of expected cost.

\begin{table}[H]\renewcommand\arraystretch{0.8}\selectfont\small
\centering
\caption{
Average percent ETC optimality gap over our 25-period test set under different demand patterns and pivoting parameters.
} \label{tab:25gaps}
\begin{tabular}{cp{1.5cm}<{\centering}p{1.5cm}<{\centering}p{1.5cm}<{\centering}p{1.5cm}<{\centering}}
\toprule
Problem Settings & ($s_t,Q_t$)-Heuristic   & ($R, S$)-Heuristic   &   ($s_t,Q$)-SDP   & ($s_t,Q$)-Heuristic\\
\midrule
\textbf{demand pattern} &&&&\\
LCY1&   2.38& 2.50& 9.56&    10.5\\
LCY2&   2.20& 2.20& 7.06&    7.60\\
SIN1&   2.52& 2.87& 6.25&    8.06\\
SIN2&   2.00& 2.03& 3.29&   3.79\\
STA&    1.45& 1.50& 1.91&     2.25\\
RAND&   2.58& 2.99& 7.24&     8.98\\
EMP1&   2.62& 3.19& 12.5&    13.3\\
EMP2&   2.50& 4.22& 24.9&    25.9\\
EMP3&   2.19& 2.79& 8.73&     9.49\\
EMP4&   2.70& 4.71& 23.8&    25.3\\
\textbf{std parameter}&&&&\\
0.1&    2.52& 2.68& 10.3&    11.4\\
0.2&    2.48& 2.50& 11.0&    11.9\\
0.3&    1.94& 3.53& 10.3&    11.3\\
\textbf{fixed ordering cost}&&&&\\
500&    2.71& 3.36& 13.8&    14.7\\
1000&   1.86& 2.61& 9.97&    10.8\\
1500&   2.35& 2.69& 7.70&     8.90\\
\textbf{penalty cost}&&&&\\
5&      2.15& 2.37& 8.79&    9.93 \\
10&     2.17& 2.97& 10.8&    11.6\\
20&     2.62& 3.37& 12.0&    13.0\\
\textbf{unit cost}&&&&\\
0&      2.53& 2.47& 11.7&    12.7\\
1&      2.10& 3.33& 9.32&     10.3\\
\midrule
\textbf{Average}&2.31&2.90&10.5&11.5\\
\bottomrule
\end{tabular}
\end{table}


\section{Conclusion}\label{sec:conclusion}

This paper investigated ($s,Q$)-type policies for the non-stationary stochastic lot-sizing problem.
By adopting a variant of \cite{bookbinder1988strategies} static-dynamic uncertainty strategy in which order quantities are fixed once and for all at the beginning of the planning horizon, we derived a stochastic dynamic formulation for the problem and proved that the associated optimal policy must take the ($s,Q$) form.

To compute optimal policy parameters, we enumerated all possible order quantity configurations to determine an optimal one, and then used a dynamic programming recursion to determine associated reorder points. Since this brute force approach is not scalable, we introduce MINLP-based heuristics to tackle large-size problems under($s,Q$)-type policies. Our heuristics leverage the MINLP approaches introduced in \citet{xiang2018computing} for the non-stationary ($s,S$) policy, in which the non-linearity of the cost function is dealt with via a piecewise linearisation of the cost function.

We carried out extensive computational experiments on a small (6-period) and a large (25-period) test set comprising 10 demand patterns and various coefficient settings.
In the numerical study on the small test set, our results show that the average optimality gaps for the ($s_t, Q_t$) policy and the ($s_t, Q$) policy versus the optimal ($s_t,S_t$)-SDP are $1.91\%$ and $3.61\%$, respectively; and the optimality gaps associated with ($s_t, Q_t$)-Heuristic and ($s_t, Q$)-Heuristic ($2.76\%$ and $4.66\%$, respectively) are close to those of the corresponding SDP.

In the numerical study on the large test set, we found that the average optimality gaps by ($s_t,Q_t$)-Heuristic remained small ($2.31\%$); while the optimality gap of the ($s_t,Q$)-Heuristic remained unsatisfactory ($11.5\%$). Our comparison against the ($R_t,S_t$)-Heuristic showed that the optimality gap of the ($s_t, Q_t$)-Heuristic was slightly better than that of the ($R_t,S_t$)-Heuristic ($2.90\%$).

Our investigation demonstrates the effectiveness of ($s, Q$)-type policies for the non-stationary stochastic lot-sizing problem. The ($s_t, Q_t$) policy can be well approximated by a heuristic and provide satisfactory results in reasonable time. The ($s_t, Q$) policy is applicable in both SDP and heuristic, while it produces larger optimality gap than ($s_t, Q_t$) policy.

\clearpage

\subfile{appendix}

\bibliography{sQ}
\clearpage

\end{document}

%% file: appendix.tex
\appendix

\linespread{1}

\renewcommand\thefigure{\Alph{section}\arabic{figure}}
\renewcommand\theequation{\Alph{section}\arabic{equation}}
\setcounter{figure}{0}
\setcounter{table}{0}
\setcounter{equation}{0}
\renewcommand\thetable{\Alph{section}\arabic{table}}

\section{Notations}\label{sec:App.notations}
\renewcommand\arraystretch{1.1}\selectfont\small
\begin{longtable}{rp{12cm}}
\caption{ Notations of important functions}\label{tab:App.notations}\\ \toprule
Functions& Explaination\\
\midrule
$c(Q)$ & cost of an order of size $Q$\\
$L_t(y)$& expected immediate holding and penalty cost when the inventory level after replenishment is $y$ at period $t$\\
$C_t(x)$ & expected total cost of the optimal policy over periods $t$ to $T$ when the opening inventory level is $x$\\
$G_t(y)$ & expected total cost over periods $t$ to $T$ when the opening inventory level is $y$ and no order is placed in period $t$\\
\midrule
$V_t(x,\vec{q}_t)$& expected total cost with a combination of reorder quantities $\vec{q}_t\in \mathcal{Q}_t$ when the opening inventory level is $x$ \\
$V_0(x)$& minimum expected total cost over $\mathcal{Q}$, the set of possible order quantities, when opening inventory level is $x$\\
$J_t(x,\vec{q}_t)$ & expected total cost with a combination of reorder quantities $\vec{q}_t\in \mathcal{Q}_t$ when no order is placed for opening inventory level $x$ in period $t$\\
$\hat{J}_t(x,\vec{q}_t)$ & expected total cost with a combination of reorder quantities $\vec{q}_t\in \mathcal{Q}_t$ when an order is placed for opening inventory level $x$ in period $t$\\
$\Delta J_t(x,\vec{q}_t)$ & $= J_t(x,\vec{q}_t)- J_t(x+Q_t,\vec{q}_t)$, the difference between expected total costs with opening inventory levels $x$ and $x+Q_t$\\
\midrule
$\mathcal{J}_t(x,\hat{\vec{q}})$ & an approximation of $J_t(x,\vec{q}_t^*)$ by MINLP\\
\bottomrule
\end{longtable}


\renewcommand\thefigure{\Alph{section}\arabic{figure}}
\renewcommand\theequation{\Alph{section}\arabic{equation}}
\setcounter{figure}{0}
\setcounter{table}{0}
\setcounter{equation}{0}
\section{MINLP model to compute $\bm{S_t}$}\label{sec:App.MINLP-Xiang}
This appendix section presents the MINLP model introduced in \citep{xiang2018computing} to compute the order-up-to level $S_t$ of the ($s_t,S_t$) policy. To properly account for the proportional ordering cost $z$, we modify the objective function in line with \citet{tarim2006modelling}. 
We apply a superscript `$S$' to distinguish decision variables from other formulations.
\vspace{-0.3cm}
\begin{alignat}{2}
\min\quad & z(\tilde{I}^S_T + \tilde{d}_{tT}) + \sum_{k=t}^T(K\delta_k^S + Q_k^S + h\cdot \tilde{H}_k + b\cdot \tilde{B}_k),& \notag \\
\mbox{s.t.}\quad
&\delta_t^S = 1,&\quad& \label{eq:sQt-MINLP-delta=1}\\
&\tilde{I}_t^S + \tilde{d}_t = S_t,&\quad& \label{eq:sQt-MINLP-delta=1-2}\\
&\delta_k^S=0 \rightarrow \tilde{I}_k^S + \tilde{d}_k = \tilde{I}_{k-1}^S, &\quad& k = t+1,\ldots, T,\label{eq:sQt-MINLP-flowBalance->delta=0}\\
&\delta_k^S=1 \rightarrow \tilde{I}_k^S + \tilde{d}_k = \tilde{I}_{k-1}^S + Q_k^S, &\quad& k = t+1,\ldots, T,\label{eq:sQt-MINLP-flowBalance->delta=1}\\
& \sum\nolimits_{j=t}^k P_{jk}^S=1, &\quad& k=t+1,\ldots,T, \label{eq:sQt-MINLP-reviewSum}\\
& P_{jk}^S \geq \delta_j^S - \sum\limits_{r=j+1}^k \delta_r^S, & \quad & k=t,\ldots,T \text{ and }j = t,\ldots, k,\label{eq:sQt-MINLP-uniqueOrder} \\
& P_{jk}^S=1 \rightarrow \tilde{H}_k = \hat{\mathcal{L}}(\tilde{I}_k^S+\tilde{d}_{jk},d_{jk}), &\quad& k=t,\ldots,T \text{ and }j = t,\ldots, k,\label{eq:sQt-MINLP-Ht-P=1}\\
& P_{jk}^S=1 \rightarrow \tilde{B}_k = \mathcal{L}(\tilde{I}_k^S+\tilde{d}_{jk},d_{jk}), &\quad& k=t,\cdots,T \text{ and }j = t,\ldots, k,\label{eq:sQt-MINLP-Bt-P=1}\\
& Q_k^S, \tilde{H}_k, \tilde{B}_k \geq 0,&\quad& k = t,\ldots, T,\label{eq:sQt-MINLP-domainHB}\\
& P_{jk}^S, \delta_k^S \in\{0,1\}, &\quad& k=t,\ldots,T\text{ and }j = t,\ldots, k.\label{eq:sQt-MINLP-domainDelta}
\end{alignat}\label{eq:modelI}
\noindent We add constraints \eqref{eq:sQt-MINLP-delta=1} and \eqref{eq:sQt-MINLP-delta=1-2} to force the system to place an order in the first period of the horizon ($t,T$) in order to approximate $S_t$. The other constraints remain as in \citep{xiang2018computing}. Constraints \eqref{eq:sQt-MINLP-flowBalance->delta=0} and \eqref{eq:sQt-MINLP-flowBalance->delta=1} capture the inventory flow balance equations and reorder conditions.
Constraint \eqref{eq:sQt-MINLP-uniqueOrder} forces $P^S_{jk} = 1$ if the most recent replenishment before period k in horizon ($t,k$) is placed in period $j$; constraint \eqref{eq:sQt-MINLP-reviewSum} ensures $P_{jk}^S = 0$ otherwise.
Constraints \eqref{eq:sQt-MINLP-Ht-P=1} and \eqref{eq:sQt-MINLP-Bt-P=1} model the expected inventory and back-ordered levels at the end of period $k$ through first order loss functions.


\renewcommand\thefigure{\Alph{section}\arabic{figure}}
\setcounter{figure}{0}
\setcounter{equation}{0}
\section{Piecewise approximation with non-stationary Poisson demand}\label{sec:App.PoissonPiecewise}
Consider a random variable $\omega$ and a scalar variable $x$, the first order loss function is defined as $\mathcal{L}(x,\omega) = \mathbb{E}[\max(\omega-x,0)]$ and its complement as $\hat{\mathcal{L}}(x,\omega) = \mathbb{E}[\max(x-\omega,0)]$. Decision variables $\tilde{H}_t\geq 0$ and $\tilde{B}_t\geq 0$ denote the expected inventory and back-order levels at the end of period $t$.

\citeauthor{rossi2014piecewise} (\citeyear{rossi2014piecewise}) presented the approach with bounding techniques to generate piecewise linear lower and upper bounds and discussed the implementation on the standard normal distribution. Instances in this paper involves non-stationary Poisson demand to enable the computation analysis on problems with small means of demand. Therefore, we extend the results of \citeauthor{rossi2014piecewise} to the Poisson distribution.

To minimise the expected inventory and back-ordere levels at the end of each period with a lower bounding piecewise linear approximation, $\tilde{H}_t$ is constrained by
\begin{equation}\label{eq:piecewise-Ht-genericDist}
\tilde{H}_t \geq (\tilde{I}_t+\sum_{j=1}^t \tilde{d}_{jt}P_{jt})\sum_{k=1}^i p_k + \sum_{j=1}^t (\sum_{k=1}^i p_k\mathbb{E}[d_{jt}|\Omega_{jt}])P_{jt},
\end{equation}
and $\tilde{B}_t$ by
\begin{equation}\label{eq:piecewise-Bt-genericDist}
\tilde{B}_t \geq -\tilde{I}_t+(\tilde{I}_t+\sum_{j=1}^t \tilde{d}_{jt}P_{jt})\sum_{k=1}^i p_k + \sum_{j=1}^t (\sum_{k=1}^i p_k\mathbb{E}[d_{jt}|\Omega_{jt}])P_{jt}.
\end{equation}
where $d_{jt}$ follows the notation in section \ref{sec:step1} denoting the convolution of $d_j$ to $d_t$, demand $d_t$ is a random variable that is of a Poisson distribution with mean $\lambda_t$, and its domain $\mathbb{R}^+$ is partitioned into $N$ disjoint adjacent subregions $\Omega_1, \Omega_2,\cdots, \Omega_N$.

According to the technique in \citep{rossi2014piecewise}, $\Omega_1 = [0, a_1]$, $\Omega_i = [a_{i-1}, a_i]$ for $i = 2,\cdots, N-1$ and $\Omega_N = [a_{N-1},\infty]$. Let the probability density function of $d_t$ be $g_{\lambda_t}(k) = e^k/ {\lambda_t !}$ and $g_{\lambda_t}^{-1}(p)$ be its inverse function, which returns the value of $k$ satisfying $g_{\lambda_t}(k) = p$, then
$$
a_i = g_{\lambda_t}^{-1}(\frac{i}{N}),
$$
and the probability $p_i$ that a realisation of the Poisson random variable $d_t$ (i.e. a value of demand $d_t$) locates within the subregion $i$ is
\begin{equation}
p_i = \text{Pr}\{d_t\in\Omega_i\} = \int_{\Omega_i}g_{\lambda_t}(u)\mathop{}\mathrm{d}u,
\end{equation}
and
\begin{equation}
\mathbb{E}[d_t|\Omega_i] = \frac{N}{i}\int_{\Omega_i} ug_{\lambda_t}(u) \mathop{} \mathrm{d}u,
\end{equation}
where $i = 1,2,\cdots,N$.

\clearpage